\documentclass{amsart}

\newtheorem{theorem}{Theorem}

\newtheorem{corollary}{Corollary} 
\newtheorem{proposition}{Proposition} 

\theoremstyle{definition}

\theoremstyle{remark}

\numberwithin{equation}{section}

\begin{document}

\title{Torsion units in integral group rings of Janko simple groups}

\author{V.A.~Bovdi}
\address{Institute of Mathematics, University of Debrecen,
P.O.  Box 12, H-4010 Debrecen, Hungary}
\email{vbovdi@math.unideb.hu}
\thanks{The research was supported by OTKA No.~K68383,
Onderzoeksraad of Vrije Universiteit Brussel,
Fonds voor Wetenschappelijk Onderzoek (Belgium),
Flemish-Polish bilateral agreement BIL2005/VUB/2006,
Francqui Stichting (Belgium) grant ADSI107
and The Royal Society of Edinburgh International Exchange Programme}

\author{E.~Jespers}
\address{Department of Mathematics, Vrije Universiteit Brussel,
Pleinlaan 2, B-1050 Brussel, Belgium}
\email{efjesper@vub.ac.be}

\author{A.B.~Konovalov}
\address{School of Computer Science, University of St Andrews,
Jack Cole Building, North Haugh, St Andrews, Fife, KY16 9SX, Scotland}
\email{alexk@mcs.st-andrews.ac.uk}


\date{April 27, 2007 and, in revised form, September 7, 2009.}

\dedicatory{Dedicated to the memory of Professor I.~S.~Luthar}


\begin{abstract}
Using the Luthar--Passi method, we investigate
the classical Zassenhaus conjecture for the
normalized unit group of integral group rings of
Janko sporadic simple groups. As a
consequence, we obtain that the  Gruenberg-Kegel
graph of the Janko groups  $J_1$, $J_2$ and $J_3$
is the same as that of the normalized unit group
of their respective integral group ring.
\end{abstract}

\maketitle

\section{Introduction  and main results}
\label{Intro}

Let $V(\mathbb Z G)$ be the normalized unit group of
the integral group ring of a finite group $G$.
A long standing conjecture \textsf{(ZC)} attributed to H.~Zassenhaus
says that every torsion unit $u\in V(\mathbb ZG)$ is
conjugate within the rational group algebra $\mathbb
Q G$ to an element in $G$ (see \cite{Zassenhaus}).

For finite simple groups, the main tool for the
investigation of the Zassenhaus conjecture is
the Luthar--Passi method, introduced in \cite{Luthar-Passi}
to solve it for $A_{5}$. Later in \cite{Hertweck1} M.~Hertweck
applied it to the investigation of \textsf{(ZC)} for $PSL(2,p^n)$
using Brauer character tables as well
as ordinary ones. In \cite{Bovdi-Grishkov-Konovalov-HeON, Bovdi-Konovalov-HS}
the method of $(p,q)$-constant characters was introduced to optimize
the Luthar--Passi method for elements of order $pq$.
The approach of Luthar and Passi proved to be
useful for groups containing non-trivial normal
subgroups as well. For some recent results on \textsf{(ZC)} we refer to
\cite{Bovdi-Hertweck}--\cite{Bovdi-Konovalov-Siciliano},
\cite{Hertweck2,Hertweck3, Hertweck1,Kimmerle}.

A weakened version of \textsf{(ZC)} can be
stated using the notion of the Gruenberg-Kegel
graph (also called the prime graph) of an arbitrary group
$X$, which is the graph $\pi (X)$ whose vertices
are labeled by primes $p$ for which there exists
an element of order $p$ in $X$ and with an edge
from $p$ to a distinct $q$ if $X$ has an element
of order $pq$. The following question is posted
by W.~Kimmerle in \cite{Oberwolfach} (Problem
21):

\centerline{ \textsf{(PQ)} \quad For a finite group
$G$, is it true that $\pi (V(\mathbb Z G)) =
\pi (G)$ ?}

\noindent Of course if \textsf{(ZC)} holds for a
finite group $G$ then \textsf{(PQ)} has a positive
answer for $G$. In \cite{Kimmerle} it is shown,
in particular, that \textsf{(PQ)} has a positive
answer for  finite Frobenius and solvable groups.
For solvable groups this result was recently 
improved by M.~Hertweck in \cite{Hertweck4},
where it was shown that orders of torsion units in
$V(\mathbb Z G)$ are exactly orders of $G$.

In \cite{Bovdi-Grishkov-Konovalov-HeON,
Bovdi-Konovalov,Bovdi-Konovalov-McL,
Bovdi-Konovalov-M23,Bovdi-Konovalov-HS,
Bovdi-Konovalov-Linton-M22,
Bovdi-Konovalov-Marcos-Suz,
Bovdi-Konovalov-Siciliano} the problem \textsf{(PQ)}
was investigated  for Mathieu, Higman-Sims,
McLaughlin, Held, Rudvalis, Suzuki and O'Nan
sporadic simple groups. In this paper we continue
these investigations for the Janko simple groups.
For $J_1$, $J_2$ and $J_3$ we again give a
positive answer to \textsf{(PQ)} and for $J_{4}$ the
question remains open. In the final section we
include some comments on the computational
difficulties encountered for $J_{4}$.

Let $G$ be a group. Put
$\mathcal{C} =\{ C_{1}, \ldots, C_{nt}, \ldots
\}$, the collection of all conjugacy classes of
$G$, where the first index denotes the order of
the elements of this conjugacy class and
$C_{1}=\{ 1\}$ (throughout the paper we will use
the ordering of conjugacy classes as used in the
GAP Character Table Library). Suppose $u=\sum
\alpha_g g \in V(\mathbb Z G)$ has finite order.
Denote by
$\nu_{nt}=\nu_{nt}(u)=\sum_{g\in
C_{nt}} \alpha_{g}$, the partial augmentation of
$u$ with respect to $C_{nt}$. From the
Berman--Higman Theorem (see, for example,
\cite{Artamonov-Bovdi}) one knows that
$\text{tr}(u) = \nu_1 =0$, and
clearly
\begin{equation}\label{E:1}
\sum_{C_{nt}\in \mathcal{C}} \nu_{nt}=1.
\end{equation}
Hence, for any character
$\chi$ of $G$, we get that $\chi(u)=\sum
\nu_{nt}\chi(h_{nt})$, where $h_{nt}$ is a
representative of a conjugacy class $ C_{nt}$.

The main results for the Janko simple groups $J_{1}$, $J_{2}$ and $J_{3}$ are as follows.

%
%

\begin{theorem}\label{T:1}

Let $G$ denote the first Janko simple
group $J_{1}$. Let $u$ be a torsion unit of
$V(\mathbb ZG)$ of order $|u|$ with
the tuple of length 15
containing partial augmentations
for all conjugacy classes of $G$.
The following properties hold.

\begin{itemize}

\item[(i)]
There are no units of orders $14$, $21$, $22$, $33$, $35$, $38$, $55$,
$57$, $77$, $95$, $133$ and $209$ in $V(\mathbb ZG)$. Equivalently,
if $|u| \not = 30$, then $|u|$ coincides with the order of some
element $g\in G$.

\item[(ii)]
If $|u| \in \{2,3,7,11,19\}$, then $u$ is rationally
conjugate to some $g\in G$.

\item[(iii)]
If $|u|=5$, then \; $\nu_{kx}=0$ \; for \; $kx \not\in\{5a,5b\}$ and
\[
\begin{split}
(\nu_{5a},\nu_{5b}) \in \{ \;  ( -1, 2),\;  ( 0, 1),\; (1, 0),\; (2,-1) \; \}.
\end{split}
\]

\item[(iv)]
 If $|u|=6$, then \; $\nu_{kx}=0$ \; for \; $kx \not\in\{2a,3a,6a\}$ and
\[
\begin{split}
(\nu_{2a},\nu_{3a},\nu_{6a}) \in \{ \;
(-4,3,2),\;(-2,0,3),\; (-2,3,0), \; (0,0,1)&,\\
(0,3,-2), \; (2,0,-1)& \; \}.
\end{split}
\]

\item[(v)]
If $|u|=10$, then \; $\nu_{kx}=0$ \; for \; $kx \not\in\{5a,5b,10a,10b\}$ and
\[
\begin{split}
(\nu_{5a},\nu_{5b},\nu_{10a},\nu_{10b}) \in \{
\; (  2, -2, 0, 1 ),\;  (  0, 0, 2, -1 ),\;
( 0, 0, 0, 1 ), \;  ( -1, 1, 1, 0 )&, \\
( 1, -1, -1, 2 ),\; (  1, -1, 1, 0 ),\;
( -2, 2, 1, 0 ),\;  ( 0, 0, -1, 2 ), \;
(  0, 0, 1, 0 )&, \\ ( -1, 1, 2, -1 ), \;
( -1, 1, 0, 1 ), \; ( 1, -1, 0, 1 ) & \; \}.
\end{split}
\]

\item[(vi)]
If $|u|=15$, then \; $\nu_{kx}=0$ \; for \; $kx \not\in\{5a,5b,15a,15b\}$ and
then
\[
\begin{split}
(\nu_{5a},\nu_{5b},\nu_{15a},\nu_{15b}) \in \{
\; (  -1, 1, 0, 1 ),\;  (  0, 0, 0, 1 ),
   (  0, 0, 1, 0 ), \; ( 1, -1, 1, 0 ) \; \}.
\end{split}
\]

\item[(vii)]
If $|u|=30$, then  \; $\nu_{kx}=0$ \; for \; $kx \not\in\{5a, 5b, 10a, 10b, 15a, 15b\}$ and
\[
\begin{split}
(\nu_{5a}, \nu_{5b}, \nu_{10a}, \nu_{10b},
\nu_{15a}, \nu_{15b}) \in \{ \; ( -1,  1, -1,
-1, 1, 2 )&, \\ ( -1,  1, -2,  0, 1, 2 ), \;
(  0,  0, -2,  0, 1, 2 ), \; (  0,  0,  0, -2,
2, 1 )&, \\ (  1, -1, -1, -1, 2, 1 ), \;
(  1, -1,  0, -2, 2, 1 ) & \; \}.
\end{split}
\]
\end{itemize}
\end{theorem}

%
%

\begin{theorem}\label{T:2}

Let $G$ denote the second Janko simple group $J_{2}$.
Let $u$ be a torsion unit of
$V(\mathbb ZG)$ of order $|u|$ with
the tuple of length 21 containing partial augmentations
for all conjugacy classes of $G$. The following properties hold.

\begin{itemize}

\item[(i)]
There are no units of orders $14$, $21$ and $35$ in $V(\mathbb ZG)$.
Equivalently, if $|u| \not \in \{ 20,24,30,40,60,120 \}$,
then $|u|$ coincides  with the order of some $g \in G$.

\item[(ii)] If $|u| \in \{7,15\}$, then $u$ is rationally conjugate to
some $g\in G$.

\item[(iii)] If $|u|=2$,
then \; $\nu_{kx}=0$ \; for \; $kx \not\in\{2a,2b\}$ and
\[
\begin{split}
(\nu_{2a},\nu_{2b}) \in
\{ \; ( 0, 1 ), \; ( -2, 3 ), \; ( 2, -1 ), \; ( 1, 0 ),  \; ( 3, -2 ), \; ( -1, 2 )  \; \}.
\end{split}
\]

\item[(iv)] If $|u|=3$,
then \; $\nu_{kx}=0$ \; for \; $kx \not\in\{3a,3b\}$ and
\[
\begin{split}
(\nu_{3a},\nu_{3b}) \in \{ \;
( 0, 1 ), \; ( 1, 0 ), \; ( -1, 2 ) \; \}.
\end{split}
\]

\item[(v)] If $|u|=4$,
then \; $\nu_{kx}=0$ \; for \; $kx \not\in\{2a,2b,4a\}$ and
\[
\begin{split}
(\nu_{2a},\nu_{2b},\nu_{4a}) \in
\{ \; ( -2, -2, 5 ), \; ( -1, -3, 5 ), \; ( -1, -1, 3 ), \; ( -1, 1, 1 )&, \\
 \; ( 0, -4, 5 ), \; ( 0, -2, 3 ), \; ( 0, 0, 1 ), \;
      ( 0, 2, -1 ), \; ( 0, 4, -3 ), \; ( 1, -3, 3 )&, \\ ( 1, -1, 1 ), \;
 ( 1, 1, -1 ), \;  ( 1, 3, -3 ), \; ( 2, 0, -1 ), \; ( 2, 2, -3 ) & \; \}.
\end{split}
\]

\item[(vi)] If $|u|=5$,
then \; $\nu_{kx}=0$ \; for \; $kx \not\in\{5a,5b,5c,5d\}$ and
\[
\begin{split}
(\nu_{5a},\nu_{5b},\nu_{5c},\nu_{5d})
\in \{ \;
(0,0,2,-1), \; (1,0,0,0), \; (0,0,0,1)&, \\
(0,0,-1,2), \; (1,0,-1,1), (0,0,1,0), \; (0,1,0,0)&, \\
(1,1,0,-1), \; (1,1,-1,0), \; (0,1,1,-1) &
\; \}.
\end{split}
\]

\item[(vii)] If $|u|=8$,
then \; $\nu_{kx}=0$ \; for \; $kx \not\in\{2a,2b,4a,8a\}$ and
\[
\begin{split}
(\nu_{2a},\nu_{2b},\nu_{4a},\nu_{8a})
\in \{ \; ( -1, -1, 0, 3 ), \; ( -1, -1, 2, 1 ), \; ( -1, 1, -2, 3 ), & \\
( -1, 1, 0, 1 ), \; ( 0, -2, 0, 3 ), \; ( 0, -2, 2, 1 ), \; ( 0, 0, -2, 3 ), \; ( 1, -1, 2, -1 )&,\\
( 0, 0, 2, -1 ), \; ( 0, 2, -2, 1 ), \; ( 0, 2, 0, -1 ), ( 0, 2, 2, -3 ), \; ( 1, -1, 0, 1 )&, \\
( 0, 0, 0, 1 ) \; ( 1, 1, -2, 1 ), \; ( 1, 1, 0, -1 ), \; ( 1, 1, 2, -3 ), \; ( 2, 0, 2, -3 ) & \; \}.
\end{split}
\]

\end{itemize}

\end{theorem}

%
%

\begin{theorem}\label{T:3}
Let $G$ denote the third Janko simple group $J_{3}$.
Let $u$ be a torsion unit of
$V(\mathbb ZG)$ of order $|u|$
with the tuple of length 21 containing partial augmentations
for all conjugacy classes of $G$. The following properties hold.

\begin{itemize}

\item[(i)]
There are no units of orders
$34$, $38$, $51$, $57$, $85$, $95$ and $323$ in $V(\mathbb ZG)$.
Equivalently, if
$|u| \not \in \{ 18, 20, 24, 30, 36, 40, 45, 60, 72, 90, 120, 180, 360 \}$,
then $|u|$ coincides with the order of some element $g \in G$.

\item[(ii)] If $|u|=2$, $u$ is rationally conjugate to
some $g\in G$.

\item[(iii)] If $|u|=3$,
then \; $\nu_{kx}=0$ \; for \; $kx \not\in\{3a,3b\}$ and
\[
\begin{split}
(\nu_{3a},\nu_{3b}) \in \{ \;
  ( 5, -4 ), \; ( 0, 1 ), \; ( -2, 3 ), \; ( 2, -1 ), \; ( -3, 4 )&, \\
  ( -4, 5 ), \; ( 1, 0 ), \; ( 3, -2 ), \; ( -1, 2 ), \; ( 4, -3 ) & \; \}.
\end{split}
\]

\item[(iv)] If $|u|=4$,
then \; $\nu_{kx}=0$ \; for \; $kx \not\in\{2a,4a\}$ and
\[
\begin{split}
(\nu_{2a},\nu_{4a}) \in \{ \;
( 0, 1 ), \; (-2,3), \; (2,-1) \; \}.
\end{split}
\]

\item[(v)] If $|u|=5$,
then \; $\nu_{kx}=0$ \; for \; $kx \not\in\{5a, 5b\}$ and
\[
\begin{split}
(\nu_{5a},\nu_{5b}) \in \{ \;
( 0, 1 ), \; ( -2, 3 ), \; ( 2, -1 ), \; ( -3, 4 ), \; ( 1, 0 )&, \\
( 3, -2 ), \; ( -1, 2 ), \; ( 4, -3 ) & \; \}.
\end{split}
\]

\item[(vi)] If $|u|=8$,
then \; $\nu_{kx}=0$ \; for \; $kx \not\in\{2a,4a,8a\}$ and
\[
\begin{split}
(\nu_{2a},\nu_{4a},\nu_{8a}) \in \{ \;
( -2, -6, 9 ), \; ( -2, -4, 7 ), \; ( -2, -2, 5 ), \; ( -2, 0, 3 )&, \\
( -2, 2, 1 ), \; ( 0, -4, 5 ), \; ( 0, -2, 3 ), \; ( 0, 0, 1 ), \;
( 0, 2, -1 ), \; ( 0, 4, -3 )&, \\
( 2, -2, 1 ), \; ( 2, 0, -1 ), \; ( 2, 2, -3 ), \; ( 2, 4, -5 ), \; ( 2, 6, -7 )& \; \}.
\end{split}
\]

\item[(vii)] If $|u|=17$,
then \; $\nu_{kx}=0$ \; for \; $kx \not\in\{17a,17b\}$ and
\[
\begin{split}
(\nu_{17a},\nu_{17b}) \in \{ \;
  ( 5, -4 ), \; ( 0, 1 ), \; ( -2, 3 ), \; ( 2, -1 ), \; ( -3, 4 )&, \\
  ( -4, 5 ), \; ( 1, 0 ), \; ( 3, -2 ), \; ( -1, 2 ), \; ( 4, -3 ) & \; \}.
\end{split}
\]

\item[(viii)] If $|u|=19$,
then \; $\nu_{kx}=0$ \; for \; $kx \not\in\{19a,19b\}$ and
\[
\begin{split}
(\nu_{19a},\nu_{19b}) \in \{ \;
  ( 5, -4 ), \; ( 0, 1 ), \; ( -2, 3 ), \; ( 2, -1 ), \; ( -3, 4 )&, \\
  ( -4, 5 ), \; ( 1, 0 ), \; ( 3, -2 ), \; ( -1, 2 ), \; ( 4, -3 ) & \; \}.
\end{split}
\]

\end{itemize}

\end{theorem}

As a consequence of the first parts of Theorems \ref{T:1} - \ref{T:3} we get the following.

%
%

\begin{corollary} If $G \in \{ J_1, J_2, J_3 \}$, then
$\pi(G)=\pi(V(\mathbb ZG))$.
\end{corollary}

\section{Preliminaries}

The following result relates the solution of
the Zassenhaus conjecture to partial augmentations
of torsion units.

\begin{proposition}\label{P:5}
(see \cite{Luthar-Passi}) Let $u\in V(\mathbb Z G)$ be
a torsion unit of order $k$. Then $u$ is
conjugate in $\mathbb QG$ to an element $g \in G$
if and only if for each $d$ dividing $k$ there is
precisely one conjugacy class $C$ with partial
augmentation $\varepsilon_{C}(u^d) \neq 0 $.
\end{proposition}

The next result already yields that several partial augmentations are zero.

\begin{proposition}\label{P:4}
(see \cite{Hertweck2}, Proposition 3.1)
Let $G$ be a finite
group and let $u$ be a torsion unit in $V(\mathbb
ZG)$. If $x$ is an element of $G$ whose $p$-part,
for some prime $p$, has order strictly greater
than the order of the $p$-part of $u$, then
$\varepsilon_x(u)=0$.
\end{proposition}

The key restriction on partial augmentations is given
by the following result.

\begin{proposition}\label{P:1}
(see \cite{Hertweck1, Luthar-Passi}) Let either
$p=0$ or $p$ be a prime divisor of $|G|$. Suppose
that $u\in V( \mathbb Z G) $ has finite order $k$
and assume $k$ and $p$ are coprime in case $p\neq
0$. If $z$ is a complex primitive $k$-th root of
unity and $\chi$ is either an ordinary character
or a $p$-Brauer character of $G$, then, for every
integer $l$, the number $$ \mu_l(u,\chi, p
)=\frac{1}{k} \sum_{d|k}Tr_{ \mathbb Q (z^d)/
\mathbb Q } \Bigl(\chi(u^d)z^{-dl}\Bigr) $$ is a
non-negative integer.
\end{proposition}

Note that if $p=0$, we will use the notation $\mu_l(u,\chi,*)$ for $\mu_l(u,\chi , 0)$.

When $s$ and $t$ are two primes such that $G$
contains no element of order $st$, and $u$ is a
normalized torsion unit of order $st$,
Proposition \ref{P:1} may be reformulated as
follows. Let $\nu_k$ be the sum of partial
augmentations of $u$ with respect to all
conjugacy classes of elements of order $k$ in
$G$, i.e. $\nu_2 = \nu_{2a}+\nu_{2b}$, etc. Then
by (\ref{E:1}) and Proposition \ref{P:4} we
obtain that $\nu_s + \nu_t = 1$ and $\nu_k = 0$
for $k \notin \{ s, t \}$. For each character
$\chi$ of $G$ (an ordinary character or a Brauer
character in characteristic not dividing $st$)
that is constant on all elements of orders $s$
and on all elements of order $t$, we have
$\chi(u) = \nu_s \chi(C_s) + \nu_t \chi(C_t)$,
where $\chi(C_r)$ denotes the value of the
character $\chi$ on any element of order $r$ from
$G$.

From the Proposition \ref{P:1} we obtain that the
values
\begin{equation}\label{E:3}
\begin{split}
\mu_l(u, \chi, p ) = \textstyle \frac{1}{st} \Bigl( \;
\chi(1) & + Tr_{\mathbb Q(z^s)/ \mathbb Q} \left(
\chi(u^s) z^{-sl} \right)  \\ & + Tr_{\mathbb
Q(z^t)/ \mathbb Q} \left( \chi(u^t) z^{-tl}
\right)
  + Tr_{\mathbb Q(z)/ \mathbb Q} \left( \chi(u) z^{-l} \right) \;
  \Bigr)
\end{split}
\end{equation}
are nonnegative integers. It follows that if
$\chi$ has the specified property, then
\begin{equation}\label{E:4}
  \mu_l(u, \chi, p) =
     \textstyle \frac{1}{st} \left( m_1 +  \nu_s m_s + \nu_t m_t \right),
\end{equation}
where
\begin{equation}\label{E:5}
\begin{split}
m_1 & = \chi(1) + \chi(C_t) Tr_{\mathbb Q(z^s)/\mathbb Q}( z^{-sl} )
              + \chi(C_s) Tr_{\mathbb Q(z^t)/\mathbb Q}( z^{-tl} ), \\
m_s & = \chi(C_s) Tr_{\mathbb Q(z)/\mathbb Q}( z^{-l} ), \qquad
m_t   = \chi(C_t) Tr_{\mathbb Q(z)/\mathbb Q}( z^{-l} ).
\end{split}
\end{equation}

Finally, we shall use well-known restrictions for torsion units.

\begin{proposition}\label{P:2}  (\cite{Cohn-Livingstone})
The order of a torsion element $u\in V(\mathbb ZG)$
divides $\text{exp}(G)$.
\end{proposition}

\begin{proposition}\label{P:6} (see \cite{Cohn-Livingstone})
Let $p$ be a prime, and let $u$ be a torsion unit of $V(\mathbb ZG)$ of order $p^n$.
Then for $m \ne n$ the sum of all partial augmentations of $u$ with respect to
conjugacy classes of elements of order $p^m$ is divisible by $p$.
\end{proposition}

\section{Proof of Theorem \ref{T:1} }

In this section we denote by $G$
the first Janko simple group $J_1$.
It is well known \cite{AFG,GAP}
that $|G|=2^3 \cdot 3 \cdot 5 \cdot 7 \cdot 11
\cdot 19$ and
$exp(G)=2 \cdot 3 \cdot 5 \cdot 7 \cdot 11 \cdot 19$.
The  character
table of $G$, as well as the Brauer character
tables for $p \in \{2,3,5,7,11,19\}$,
denoted by $\mathfrak{BCT}{(p)}$,
can be found using the
computational algebra system GAP \cite{GAP}, which
derives its data from \cite{AFG,ABC}.
We will use the
notation, including the indices, for the
characters and conjugacy classes as used in the
GAP Character Table Library.

Since $G$ only possesses elements of orders
$2$, $3$, $5$, $6$, $7$, $10$, $11$, $15$ and
$19$, we first investigate normalized units of
these orders. After this, by Proposition \ref{P:2}, the
order of each torsion unit divides $exp(G)$,
so it is enough to consider normalized units of orders
$14$, $21$, $22$, $30$, $33$, $35$, $38$, $55$,
$57$, $77$, $95$, $133$ and $209$,
because if $u$ is a unit of another possible order, then there is
$t \in \mathbb N$ such that $u^t$ has an order from this list.
We shall prove
that units of all these orders except $30$
do not appear in $V(\mathbb ZG)$.

Assume that $u$ is a non-trivial normalized unit and consider each case separately.

\noindent $\bullet$ Let $|u| \in \{2,3,7,11\}$.
Since there is only one conjugacy class in $G$
consisting of elements of order $|u|$, this
case follows at once from Propositions \ref{P:5} and \ref{P:4}.

\noindent$\bullet$ Let $|u|=19$.
By (\ref{E:1}) and Proposition \ref{P:4} we have that
$\nu_{19a}+\nu_{19b}+\nu_{19c}=1$.
Applying
Proposition \ref{P:1} to characters $\chi_{2},\chi_{8},\chi_{13}$ in
$\mathfrak{BCT}{(11)}$ we get the system:

\[
\begin{split}
\mu_1(u,\chi_{2},11) & = \textstyle\frac{1}{19} ( -t_1 +   7) \geq 0; \quad 
\mu_1(u,\chi_{8},11)   = \textstyle\frac{1}{19} ( t_1 +  69) \geq 0; \\ 
\mu_2(u,\chi_{2},11) & = \textstyle\frac{1}{19} ( -t_2 +   7) \geq 0; \quad 
\mu_2(u,\chi_{8},11)   = \textstyle\frac{1}{19} ( t_2 +  69) \geq 0; \\ 
\mu_4(u,\chi_{2},11) & = \textstyle\frac{1}{19} (t_3 +   7) \geq 0; \qquad 
\mu_4(u,\chi_{8},11)   = \textstyle\frac{1}{19} (-t_3 +  69) \geq 0; \\ 
\mu_1(u,&\chi_{13},11)  = \textstyle\frac{1}{19} ( 14 \nu_{19a} -  5 \nu_{19b} -  5 \nu_{19c} + 119) \geq 0; \\ 
\mu_2(u,&\chi_{13},11)  = \textstyle\frac{1}{19} ( -5 \nu_{19a} + 14 \nu_{19b} -  5 \nu_{19c} + 119) \geq 0; \\ 
\mu_4(u,&\chi_{13},11)  = \textstyle\frac{1}{19} ( -5 \nu_{19a} -  5 \nu_{19b} + 14 \nu_{19c} + 119) \geq 0. \\ 
\end{split}
\]
where
$t_1= 7 \nu_{19a} - 12 \nu_{19b} +  7 \nu_{19c}$,
$t_2= 7 \nu_{19a} +  7 \nu_{19b} - 12 \nu_{19c}$ and
$t_3= 12 \nu_{19a} -  7 \nu_{19b} -  7 \nu_{19c}$.
From these restrictions and the requirement that all
$\mu_i(u,\chi_{j},p)$ must be non-negative integers
we get that
$(\nu_{19a},\nu_{19b},\nu_{19c}) \in
\{(1,0,0),(0,1,0),(0,0,1)\}$.

Thus, for units of orders 2, 3, 7, 11 and 19
there is precisely one conjugacy
class with non-zero partial augmentation
so Proposition \ref{P:5} yields part (ii) of Theorem \ref{T:1}.

Note that using the LAGUNA package \cite{LAGUNA} in combination with
constraint solvers MINION \cite{MINION} and ECLiPSe \cite{ECLiPSe},
we computed inequalities from Proposition \ref{P:1} for every
irreducible character from ordinary and Brauer character tables,
and for every $0 \leq l \leq |u|-1$ (it is enough to enumerate $l$
in this range since $z^{|u|}=1$, so for bigger values of $l$ we
will not have new inequalities), but the only inequalities that really
matter are those listed above. The same remark applies for
all other orders of torsion units considered in the paper.

\noindent $\bullet$ Let $u$ be a unit of order
$5$. By (\ref{E:1}) and Proposition \ref{P:4} we get
$\nu_{5a}+\nu_{5b}=1$.
Again applying Proposition \ref{P:1} to characters in
$\mathfrak{BCT}{(11)}$ we get the following system of inequalities:
\[
\small{
\begin{split}
\mu_1(u,\chi_2,11)   &= \textstyle \frac{1}{5} ( 3 \nu_{5a} - 2 \nu_{5b} + 7   ) \geq 0; \quad \;
\mu_2(u,\chi_2,11)    = \textstyle \frac{1}{5} (-2 \nu_{5a} + 3 \nu_{5b} + 7  ) \geq 0; \\
\mu_1(u,\chi_3,11)   &= \textstyle \frac{1}{5} (-4 \nu_{5a} +   \nu_{5b} + 14 ) \geq 0; \; \;
\mu_2(u,\chi_3,11)    = \textstyle \frac{1}{5} (   \nu_{5a} - 4 \nu_{5b} + 14 ) \geq 0; \\
\mu_1(u,\chi_5,11)   &= \textstyle \frac{1}{5} (   \nu_{5a} - 4 \nu_{5b} + 49 ) \geq 0; \quad \;
\mu_1(u,\chi_6,11)    = \textstyle \frac{1}{5} (-6 \nu_{5a} + 4 \nu_{5b} + 56 ) \geq 0; \\
\mu_2(u,\chi_6,11)   &= \textstyle \frac{1}{5} ( 4 \nu_{5a} - 6 \nu_{5b} + 56 ) \geq 0; \quad
\mu_2(u,\chi_7,11)    = \textstyle \frac{1}{5} (-4 \nu_{5a} + 6 \nu_{5b} + 64 ) \geq 0; \\
\mu_2(u,\chi_8,11)   &= \textstyle \frac{1}{5} (   \nu_{5a} - 4 \nu_{5b} + 69 ) \geq 0; \quad \;
\mu_1(u,\chi_{12},11) = \textstyle \frac{1}{5} ( 4 \nu_{5a} -   \nu_{5b} + 106) \geq 0, \\
\end{split}}
\]
that has only four integer solutions
$(\nu_{5a},\nu_{5b}) \in \{
(0,1),(2,-1),(1,0),(-1,2) \}$ such that all
$\mu_i(u,\chi_{j},11)$ are non-negative integers,
so part (iii) of Theorem \ref{T:1} is proved.

\noindent$\bullet$ Let $|u|=6$.
By (\ref{E:1}) and Proposition \ref{P:4} we have that
$\nu_{2a}+\nu_{3a}+\nu_{6a}=1$.
Applying Proposition \ref{P:1} to characters in
$\mathfrak{BCT}{(11)}$ we get the following system of
inequalities:

\[
\begin{split}
\mu_3(u,\chi_4,11) &= \textstyle
\frac{1}{6} ( -6 \nu_{2a} + 24 ) \geq 0; \quad \,
\mu_0(u,\chi_4,11) = \textstyle \frac{1}{6} ( 6
\nu_{2a} + 30 ) \geq 0; \\
\mu_0(u,\chi_6,11) &=
\textstyle \frac{1}{6} ( 4 \nu_{3a} + 60 ) \geq
0; \qquad \mu_3(u,\chi_6,11) = \textstyle
\frac{1}{6} (-4 \nu_{3a} + 60 ) \geq 0; \\
\mu_0(u,\chi_2,11) &= \textstyle \frac{1}{6} ( -t_1 +8 ) \geq 0; \qquad \quad
\mu_3(u,\chi_2,11)  = \textstyle \frac{1}{6} ( t_1 + 10 ) \geq
0; \\
\mu_0(u,\chi_3,11) &= \textstyle \frac{1}{6} ( -2 t_2 + 10 ) \geq 0; \quad \quad
\mu_1(u,\chi_3,11) = \textstyle \frac{1}{6} ( -t_2 + 17 ) \geq 0; \\
& \qquad \mu_3(u,\chi_3,11) = \textstyle \frac{1}{6} ( 2 t_2 + 14 ) \geq 0, \\
\end{split}
\]
where $t_1=2 \nu_{2a} -2 \nu_{3a} + 2 \nu_{6a}$ and $t_2=2 \nu_{2a} + \nu_{3a}- \nu_{6a}$
which has only six integer solutions such that all
$\mu_i(u,\chi_{j},11)$ are non-negative integers.
These are as  listed in part (iv) of Theorem \ref{T:1}.

\noindent$\bullet$ Let $|u|=10$.
By (\ref{E:1}) and Proposition \ref{P:4} we get
$\nu_{2a}+\nu_{5a}+\nu_{5b}+\nu_{10a}+\nu_{10b}=1$.
For any character $\chi$ of $G$
we need to consider 4 cases
from part (iii) of Theorem \ref{T:1}:
\[
\begin{matrix}
\text{Case 1.} & \chi(u^2) & = & \chi(5a). \qquad &
\text{Case 3.} & \chi(u^2) & = & 2 \chi(5a) - \chi(5b). \\
\text{Case 2.} & \chi(u^2) & = & \chi(5b). \qquad &
\text{Case 4.} & \chi(u^2) & = & - \chi(5a) + 2 \chi(5b).
\end{matrix}
\]
Here and below, $\chi(5a)$ denotes the
value of the character $\chi$ on the
representative of the conjugacy class $C_{5a}$,
etc.

Applying Proposition \ref{P:1} to characters in
$\mathfrak{BCT}{(11)}$ we get the system:
\[
\begin{split}
\mu_0(u,\chi_{2},11) & = \textstyle \frac{1}{10}(-t_1 + 4) \geq 0; \quad \quad
\mu_5(u,\chi_{2},11)   = \textstyle \frac{1}{10}(t_1 + 6) \geq 0; \\
\mu_2(u,\chi_{3},11) & = \textstyle \frac{1}{10}(t_2 + \alpha_1) \geq 0; \quad \quad
\mu_3(u,\chi_{3},11)   = \textstyle \frac{1}{10}(-t_2 + \alpha_2) \geq 0; \\
\mu_0(u,\chi_{6},11) & = \textstyle \frac{1}{10}(t_3 + 60) \geq 0; \qquad
\mu_5(u,\chi_{6},11)   = \textstyle \frac{1}{10}(-t_3 + 60) \geq 0; \\
\mu_1(u,\chi_{6},11) & = \textstyle \frac{1}{10}(-t_4 + \alpha_3) \geq 0; \quad
\mu_4(u,\chi_{6},11)   = \textstyle \frac{1}{10}(t_4 + \alpha_3) \geq 0, \\
\end{split}
\]
where
$t_1 = 4 \nu_{2a} + 2 \nu_{5a} + 2 \nu_{5b} - 6 \nu_{10a} - 6 \nu_{10b}$, \;
$t_2 = 2 \nu_{2a} - 4 \nu_{5a} +  \nu_{5b} + 2 \nu_{10a} - 3 \nu_{10b}$,
$t_3 = 4 \nu_{5a} + 4 \nu_{5b}$ and
$t_4 = 4 \nu_{5a} - 6 \nu_{5b}$,
and $(\alpha_1,\alpha_2,\alpha_3)$ is equal
to $(13,17,50)$, $(8,12,60)$, $(18,22,40)$ and
$(3,7,70)$ in cases 1-4 respectively.

Now denote $t_1 = 2 \nu_{2a} + \nu_{5a} +
\nu_{5b} - 3 \nu_{10a} - 3 \nu_{10b} \in \mathbb Z$, then from
the first two inequalities (recall that all $\mu_i(u,\chi_{j},p)$
must be non-negative integers) we obtain that $t_1
\in \{ -3, 2 \}$. Put $t_2 = 2 \nu_{2a} - 4
\nu_{5a} +  \nu_{5b} + 2 \nu_{10a} - 3
\nu_{10b}  \in \mathbb Z$. From the third and fourth
inequalities $t_2$ belongs to the set $\{ -13,
-3, 7, 17 \}$, $\{ -8, 2, 12 \}$, $\{ -18, -8, 2,
12, 22 \}$  and $\{ -3, 7 \}$ in cases 1-4
respectively. Put $t_3 = \nu_{5a}+\nu_{5b} \in \mathbb Z$.
From the fifth and eighth inequalities it follows
that $t_3 \in \{ 5k \mid -3 \le k \le 3 \}$.

Finally, put $t_4 = 2 \nu_{5a} - 3 \nu_{5b}$.
Considering the sixth and seventh inequalities we get
that $t_4$ belongs to the set $\{ 5k \mid -5 \le k \le 5 \}$, $\{ 5k
\mid -6 \le k \le 6 \}$, $\{ 5k \mid -4 \le k \le
4 \}$ and $\{ 5k \mid -7 \le k \le 7 \}$ in cases
1-4 respectively.

We obtain the system of five linear equations:
\[
\begin{split}
\nu_{2a}+\nu_{5a}+\nu_{5b}+\nu_{10a}+\nu_{10b} & = 1 ;\\
2 \nu_{2a} + \nu_{5a} + \nu_{5b} - 3 \nu_{10a} - 3 \nu_{10b} & = t_1 ;\\
2 \nu_{2a} - 4 \nu_{5a} +  \nu_{5b} + 2 \nu_{10a} - 3 \nu_{10b} & = t_2 ;\\
\nu_{5a}+\nu_{5b} = t_3;\qquad 2 \nu_{5a} - 3 \nu_{5b} & = t_4 .\\
\end{split}
\]
Since the matrix of the system is
non-degenerate, it has a unique
solution for any values of parameters $t_i$.
For each of the allowable values of $t_1$,
$t_2$, $t_3$ and $t_4$ we can thus compute
the unique integer solution of this system
of equations.

Now Proposition \ref{P:1} for
$\mathfrak{BCT}{(11)}$ also gives the following
additional inequalities:

\[
\begin{split}
\mu_2(u,\chi_{2},11) & = \textstyle \frac{1}{10} ( t_1 + \beta_2) \geq 0; \quad   
\mu_3(u,\chi_{2},11)   = \textstyle \frac{1}{10} (-t_1 + \beta_3) \geq 0; \\   
\mu_1(u,\chi_{2},11) & = \textstyle \frac{1}{10} (-t_2 + \beta_1) \geq 0; \quad   
\mu_4(u,\chi_{2},11)   = \textstyle \frac{1}{10} ( t_2 + \beta_4) \geq 0; \\   
\mu_1(u,\chi_{3},11) & = \textstyle \frac{1}{10} (-t_3 + \beta_5) \geq 0; \quad 
\mu_4(u,\chi_{3},11)   = \textstyle \frac{1}{10} (t_3+ \beta_6) \geq 0, \\  
\mu_0(u,\chi_{3},11) & = \textstyle \frac{1}{10} (-8 \nu_{2a} + 6 \nu_{5a} + 6 \nu_{5b} + 2 \nu_{10a} + 2 \nu_{10b} + 18) \geq 0; \\     
\end{split}
\]
$t_1=\nu_{2a} + 3 \nu_{5a} - 2 \nu_{5b} +  \nu_{10a} - 4 \nu_{10b}$,
$t_2=\nu_{2a} - 2 \nu_{5a} + 3 \nu_{5b} - 4 \nu_{10a} +  \nu_{10b}$
and
$t_3=2 \nu_{2a} +  \nu_{5a} - 4 \nu_{5b} - 3 \nu_{10a} + 2 \nu_{10b}$,
and \quad
$(\beta_1,\beta_2,\beta_3,\beta_4,\beta_5,\beta_6)$ \quad is equal to \quad
$(11,-4,6,9,12,8)$, \; $(6,9,11,4,17,13)$,
$(16,-1,1,14,7,3)$, \; $(1,14,16,-1,22,18)$ \quad
in cases 1-4 respectively.
It then follows that $\nu_{2a}=0$ and the only integer solutions with
non-negative integers $\mu_i(u,\chi_{j},11)$ are those
listed in part (v) of Theorem \ref{T:1}.

\noindent$\bullet$ Let $|u|=15$.
By (\ref{E:1}) and Proposition \ref{P:4} we get
$\nu_{3a}+\nu_{5a}+\nu_{5b}+\nu_{15a}+\nu_{15b}=1$.
Since $|u^3|=5$, for any character $\chi$ of $G$
we need to consider four cases,
defined by part (iii) of Theorem \ref{T:1}:
\[
\begin{matrix}
\text{Case 1.} & \chi(u^3) & = &
\chi(5a). \qquad & \text{Case 3.} & \chi(u^3) & =
& 2 \chi(5a) -   \chi(5b). \\ \text{Case 2.} &
\chi(u^3) & = & \chi(5b). \qquad &\text{Case 4.}
& \chi(u^3) & = & - \chi(5a) + 2 \chi(5b).
\end{matrix}
\]
Again applying Proposition \ref{P:1} to characters in
$\mathfrak{BCT}{(11)}$ we get the system:
\[
\begin{split}
\mu_0(u,\chi_{2},11) & = \textstyle \frac{1}{15}(-2 t_1 + 7) \geq 0; \quad
\mu_5(u,\chi_{2},11)   = \textstyle \frac{1}{15}(t_1 + 4) \geq 0; \\
\mu_0(u,\chi_{3},11) & = \textstyle \frac{1}{15}(-2 t_2 + 18) \geq 0; \quad
\mu_5(u,\chi_{3},11)   = \textstyle \frac{1}{15}(t_2 + 21) \geq 0; \\
\mu_1(u,\chi_{4},11) & = \textstyle \frac{1}{15}(t_3 + \alpha_1) \geq 0; \quad
\mu_6(u,\chi_{4},11)   = \textstyle \frac{1}{15}(-2 t_3 + \alpha_1) \geq 0; \\
\mu_1(u,\chi_{6},11) & = \textstyle \frac{1}{15}(t_4 + \alpha_2) \geq 0; \quad
\mu_1(u,\chi_{7},11)   = \textstyle \frac{1}{15}(-t_4 + \alpha_3) \geq 0, \\
\end{split}
\]
where
$t_1=-4\nu_{3a} + 2 \nu_{5a} + 2 \nu_{5b} - 4 \nu_{15a} - 4 \nu_{15b}$,
$t_2=4\nu_{3a} - 6 \nu_{5a} - 6 \nu_{5b} - 6 \nu_{15a} - 6 \nu_{15b}$,
$t_3=- 3 \nu_{5a} + 2 \nu_{5b} + 5 \nu_{15b}$ and
$t_4=2\nu_{3a} - 4 \nu_{5a} + 6 \nu_{5b} + 2 \nu_{15a} - 3 \nu_{15b}$,
and $(\alpha_1,\alpha_2,\alpha_3)$ is equal
to $(25,48,72)$,$(30,58,62)$,$(20,38,82)$ and
$(35,68,52)$ in cases 1-4 respectively.

Now denote $t_1 = 2 \nu_{3a} - \nu_{5a} -
\nu_{5b} + 2 \nu_{15a} + 2 \nu_{15b} \in \mathbb Z$, then from
the first two inequalities (again recall that all $\mu_i(u,\chi_{j},p)$
must be non-negative integers)
we obtain that
$t_1=2$. Put
$t_2 = 2 \nu_{3a} - 3 \nu_{5a} - 3
\nu_{5b} - 3 \nu_{15a} - 3 \nu_{15b} \in \mathbb Z$. From
the third and fourth inequalities $t_2=-3$. Put
$t_3 = 3 \nu_{5a} - 2 \nu_{5b} - 5 \nu_{15b}\in\mathbb Z$.
From the fifth and sixth inequalities it follows
that $t_3$ belongs to the set $\{-5,10,25\}$,
$\{-15,0,15,30\}$, $\{-10,5,20\}$ and
$\{-10,5,20,35\}$ in cases 1-4 respectively.

Finally, put $t_4 = 2\nu_{3a} - 4 \nu_{5a} + 6
\nu_{5b} + 2 \nu_{15a} - 3 \nu_{15b} \in \mathbb
Z$. Considering the seventh and eighth
inequalities we get that $t_{4}$ belongs to the
set $\{\; \gamma+15k \; \mid \; k=0,\dots,8 \;
\}$, where $\gamma$ is equal to $-48$, $-58$,
$-38$, $-68$ in cases 1-4 respectively. So
we get
\[
\begin{split}
\nu_{3a}+\nu_{5a}+\nu_{5b}+\nu_{15a}+\nu_{15b} &
= 1 ;\\ 2 \nu_{3a} - \nu_{5a} - \nu_{5b} + 2
\nu_{15a} + 2 \nu_{15b} & = 2 ;\\ 2 \nu_{3a} - 3
\nu_{5a} - 3 \nu_{5b} - 3 \nu_{15a} - 3 \nu_{15b}
& = -3 ;\\ 3 \nu_{5a} - 2 \nu_{5b} - 5 \nu_{15b}
& = t_3 ;\\ 2\nu_{3a} - 4 \nu_{5a} + 6 \nu_{5b} +
2 \nu_{15a} - 3 \nu_{15b} & = t_4 .\\
\end{split}
\]
Since the matrix of the system is
non-degenerate, again such a system has a unique
solution for any values of the parameters $t_i$. For each
of the allowable values of $t_3$ and $t_4$, we thus can
compute the unique integer solution of this system of
equations.

Now Proposition \ref{P:1} for $\mathfrak{BCT}{(11)}$ also
gives the following additional inequalities:
\[
\begin{split}
\mu_2(u,\chi_{2},11) & = \textstyle \frac{1}{15} (t_1 + \beta_2) \geq 0; \quad 
\mu_3(u,\chi_{2},11)   = \textstyle \frac{1}{15} (-t_1 + \beta_3) \geq 0; \\ 
\mu_1(u,\chi_{2},11) & = \textstyle \frac{1}{15} (t_2 + \beta_1) \geq 0; \quad 
\mu_6(u,\chi_{2},11)   = \textstyle \frac{1}{15} (-2 t_2 + \beta_4) \geq 0; \\ 
\mu_1(u,\chi_{3},11) & = \textstyle \frac{1}{15} (t_3 + \beta_5) \geq 0; \quad 
\mu_6(u,\chi_{3},11)   = \textstyle \frac{1}{15} (-2 t_3 + \beta_7) \geq 0; \\ 
\mu_3(u,\chi_{3},11) & = \textstyle \frac{1}{15} (2\nu_{3a} - 8 \nu_{5a} + 2 \nu_{5b} - 8 \nu_{15a} + 2 \nu_{15b} + \beta_6) \geq 0, \\ 
\end{split}
\]
where
$t_1=\nu_{3a} - 3 \nu_{5a} + 2 \nu_{5b} + 6 \nu_{15a} - 4 \nu_{15b}$,
$t_2=\nu_{3a} + 2 \nu_{5a} - 3 \nu_{5b} - 4 \nu_{15a} + 6 \nu_{15b}$
and
$t_3=-\nu_{3a} -  \nu_{5a} + 4 \nu_{5b} -  \nu_{15a} + 4 \nu_{15b} $,
and
$(\beta_1,\beta_2,\beta_3,\beta_4,\beta_5,\beta_6,\beta_7)$
is equal to $(9,4,7,12,11,13,8)$,
$(4,9,12,7,16,8,13)$, \; $(14,-1,2,17,6,18,3)$ \; and
$(-1,14,17,2,21,3,18)$ \; in cases 1-4 respectively.
It then follows that $\nu_{3a}=0$ and the only integer solutions with
non-negative integers $\mu_i(u,\chi_{j},11)$ are those
listed in part (vi) of Theorem \ref{T:1}.

$\bullet$ Let $|u|=30$. Since $|u^i|=\frac{|u|}{(|u|,i)}$,
by parts (iii)-(vi) of Theorem \ref{T:1} there are 4, 6, 12 and 4
tuples of partial augmentations for orders  5, 6, 10 and 15 respectively,
so we need to consider $4 \cdot 6 \cdot 12 \cdot 4 = 1152$ cases. Using the
LAGUNA package \cite{LAGUNA} together with MINION
and ECLiPSe \cite{ECLiPSe, MINION}, we
constructed and solved all of them, and
only six cases given in the table
below yield a non-trivial solution (see
part (vii) of Theorem~\ref{T:1}).
$$
\small{
\begin{array}{|c|c|c|c|}\hline
\chi(u^6) & \chi(u^5)   & \chi(u^3) & \chi(u^2) \\ \hline
\chi(5b)           & 3\chi(3a)-2\chi(2a) & 2\chi(5b)-2\chi(5a)+\chi(10a) & \chi(15a) \\ 
\chi(5b)           & 3\chi(3a)-2\chi(2a) & 2\chi(10b)-\chi(10a)          & \chi(15a) \\ 
\chi(5a)           & 3\chi(3a)-2\chi(2a) &  2\chi(10a)-\chi(10b)          & \chi(15b) \\ 
\chi(5a)           & 3\chi(3a)-2\chi(2a) &  2\chi(5a)-2\chi(5b)+\chi(10b) & \chi(15b) \\ 
2\chi(5a)-\chi(5b) & 3\chi(3a)-2\chi(2a) & \chi(10a)                     & \chi(5a)-\chi(5b)+\chi(15a) \\ 
2\chi(5b)-\chi(5a) & 3\chi(3a)-2\chi(2a) & \chi(10b)                     & \chi(5b)-\chi(5a)+\chi(15b) \\ \hline 
\end{array}}
$$

$\bullet$ It remains to prove part (i) of Theorem \ref{T:1}, that is to
show that $V(\mathbb ZG)$ has no elements of orders
$14$, $21$, $22$, $33$, $35$, $38$, $55$, $57$,
$77$, $95$, $133$ and $209$.
We give a detailed proof for order 33. Other cases can be
derived similarly from the table below containing the data for
the constraints on partial augmentations $\nu_p$ and $\nu_q$
for possible orders $pq$ (including order 33 as well)
accordingly to (\ref{E:3})--(\ref{E:5}).

If $|u|=33$, then $\nu_{3}+\nu_{11}=1$.
Consider ordinary characters $\xi=\chi_{7}$ and $\tau=\chi_{6}$,
which are encoded in the table as $\xi=(7)_{[*]}$ and $\tau=(6)_{[*]}$
respectively. These characters
are constant on elements of order 3 and elements of order 11:
$\xi(C_3) = 2$, $\xi(C_{11}) = 0$, $\tau(C_3) = -1$ and $\tau(C_{11}) = 0$.
Now we obtain the following system:
\[
\begin{split}
\mu_{0}(u,\xi,*)   =   \textstyle \frac{1}{33} & (40 \nu_{3} + 81) \geq 0; \qquad
\mu_{1}(u,\xi,*)   =   \textstyle \frac{1}{33}(2 \nu_{3} + 75) \geq 0; \\
& \mu_{0}(u,\tau,*) =  \textstyle \frac{1}{33}(-20 \nu_{3} + 75) \geq 0, \\
\end{split}
\]
which has no integral solution $(\nu_{3},\nu_{11})$
such that all $\mu_i(u,\chi_{j},*)$ are non-negative integers.

To complete the proof, we give
the data for part (i) of Theorem \ref{T:1} in the table below
(we use the notation $\xi=(i)_{[p]}$ for $p$-Brauer characters).

\centerline{\small{
\begin{tabular}{|c|c|c|c|c|c|c|c|c|c|c|c|c|c|}
\hline
$|u|$&$p$&$q$&$\xi, \; \tau$&$\xi(C_p)$&$\xi(C_q)$&$l$&$m_1$&$m_p$&$m_q$ \\
\hline
   &   &   &                       &    &   & 0 & 82   & 30   & 0 \\
14 & 2 & 7 & $\xi=(6)_{[*]}$       & 5  & 0 & 1 & 72   & 5    & 0 \\
   &   &   &                       &    &   & 7 & 72   & -30  & 0 \\
\hline
   &   &   &                       &   &   & 0 & 60   & 24    & 0 \\
21 & 3 & 7 & $\xi=(2)_{[*]}$       & 2 & 0 & 1 & 54   & 2     & 0 \\
   &   &   &                       &   &   & 7 & 54   & -12   & 0 \\
\hline
   &   &   &                       &   &   & 0 & 82   & 50    & 0 \\
22 & 2 & 11& $\xi=(6)_{[*]}$       & 5 & 0 & 1 & 72   & 5     & 0 \\
   &   &   &                       &   &   & 7 & 72   & -50   & 0 \\
\hline
   &   &   & $\xi=(7)_{[*]}$       & 2 & 0 & 0 & 81   & 40    & 0 \\
33 & 3 & 11& $\xi=(7)_{[*]}$       & 2 & 0 & 1 & 75   & 2     & 0 \\
   &   &   & $\tau=(6)_{[*]}$      & -1& 0 & 0 & 75   & -20   & 0 \\
\hline
   &   &   & $\xi=(6)_{[*]}$       & 2 & 0 & 0 & 85   & 48    & 0 \\
35 & 5 & 7 & $\xi=(6)_{[*]}$       & 2 & 0 & 1 & 75   & 2     & 0 \\
   &   &   & $\tau =(12)_{[*]}$    & -2& 0 & 0 & 125  & -48   & 0 \\
\hline
   &   &   & $\xi=(15)_{[*]}$      & 1 & 0 & 1 & 208  & 1     & 0 \\
38 & 2 & 19& $\tau =(4)_{[*]}$     & 4 & 0 & 0 & 80   & 72    & 0 \\
   &   &   & $\tau =(4)_{[*]}$     & 4 & 0 & 19& 72   & -72   & 0 \\
\hline
   &   &   & $\xi=(6)_{[*]}$       & 2 & 0  & 11& 75   & -20   & 0 \\
55 & 5 & 11& $\tau =(9)_{[*]}$     & 0 & -1 & 0 & 110  & 0     & -40 \\
   &   &   & $\tau =(9)_{[*]}$     & 0 & -1 & 5 & 121  & 0     & 4 \\
\hline
   &   &   &                       &   &   & 0 & 72   & -72    & 0 \\
57 & 3 & 19& $\xi=(7)_{[2]}$       & -2& 0 & 1 & 78   & -2     & 0 \\
   &   &   &                       &   &   & 19& 78   & 36     & 0 \\
\hline
   &   &   &                       &   &   & 0 & 66   & 0     & 60 \\
77 & 7 & 11& $\xi=(3)_{[2]}$       & 0 & 1 & 7 & 55   & 0     & -6 \\
   &   &   &                       &   &   & 11& 66   & 0     & -10 \\
\hline
   &   &   &                       &   &   & 0 & 38   & 0     & 72 \\
95 & 5 & 19& $\xi=(2)_{[2]}$       & 0 & 1 & 5 & 19   & 0     & -4 \\
   &   &   &                       &   &   & 19& 38   & 0     & -18 \\
\hline
   &   &   &                       &   &   & 0 & 38   & 0     & -108 \\
133& 7 & 19& $\xi=(3)_{[2]}$       & 0 & -1& 7 & 57   & 0     & 6 \\
   &   &   &                       &   &   & 19& 38   & 0     & 18 \\
\hline
   &   &   &                       &   &   & 0 & 66   & -180  & 0 \\
209& 11& 19& $\xi=(7)_{[2]}$       & -1& 0 & 11& 66   & 10    & 0 \\
   &   &   &                       &   &   & 19& 77   & 18    & 0 \\
\hline
\end{tabular}
}}

\section{Proof of Theorem \ref{T:2}}

Let $G$ be the second Janko simple group $J_2$.
It is well known \cite{AFG,GAP} that
$|G|=2^7 \cdot 3^3 \cdot 5^2 \cdot 7$ and
$exp(G)=2^3 \cdot 3 \cdot 5 \cdot 7$.

Since the group $G$ only possesses elements of orders
$2$, $3$, $4$, $5$, $6$, $7$, $8$, $10$, $12$ and $15$,
we will first investigate normalized units of
these orders. Due to Proposition \ref{P:2}, the
order of each torsion unit divides the exponent
of $G$, so it remains to consider normalized units of orders
$14$, $20$, $21$, $24$, $30$ and $35$.
We shall prove that units of all these orders except
$20$, $24$ and $30$
do not appear in $V(\mathbb ZG)$.

Assume that $u$ is a non-trivial normalized unit and consider each case separately.

\noindent $\bullet$ Let $|u|=2$. By
(\ref{E:1}) and Proposition \ref{P:4} we get
$\nu_{2a}+\nu_{2b}=1$. By Proposition \ref{P:1}
we get the system of inequalities:
\[
\begin{split}
\mu_{0}(u,\chi_{2},*) & = \textstyle \frac{1}{2} (-2 \nu_{2a} + 2 \nu_{2b} + 14) \geq 0; \quad
\mu_{1}(u,\chi_{2},*)   = \textstyle \frac{1}{2} (2 \nu_{2a} - 2 \nu_{2b} + 14) \geq 0; \\
\mu_{0}(u,\chi_{4},*) & = \textstyle \frac{1}{2} (5 \nu_{2a} - 3 \nu_{2b} + 21) \geq 0; \quad
\mu_{1}(u,\chi_{4},*)   = \textstyle \frac{1}{2} (-5 \nu_{2a} + 3 \nu_{2b} + 21) \geq 0, \\
\end{split}
\]
that has only six integer solutions $(\nu_{2a},\nu_{2b})$
as listed in part (iii) of Theorem \ref{T:2}.

\noindent $\bullet$ Let $|u|=3$. By
(\ref{E:1}) and Proposition \ref{P:4} we get
$\nu_{2a}+\nu_{2b}=1$. Again using Proposition
\ref{P:1} we get the system of inequalities:
\[
\begin{split}
\mu_{0}(u,\chi_{2},*) & = \textstyle \frac{1}{3} (10 \nu_{3a} - 2 \nu_{3b} + 14) \geq 0; \\
\mu_{1}(u,\chi_{2},*) & = \textstyle \frac{1}{3} (-5 \nu_{3a} +  \nu_{3b} + 14) \geq 0; \\
\mu_{0}(u,\chi_{2},2) & = \textstyle \frac{1}{3} (-6 \nu_{3a} + 6) \geq 0, \\
\end{split}
\]
that has only three integer solutions $(\nu_{3a},\nu_{3b})$ as
listed in part (iv) of Theorem \ref{T:2}.

\noindent $\bullet$ Let $|u|=4$. By
(\ref{E:1}) and Proposition \ref{P:4} we get $
\nu_{2a}+\nu_{2b}+\nu_{4a}=1$. We need to
consider six cases defined by part (iii) of
Theorem \ref{T:2}. In each of these cases, we
apply Proposition \ref{P:1} to get the following
systems of inequalities:
\[
\begin{split}
\mu_{0}(u,\chi_{2},*) & = \textstyle \frac{1}{4} (-4 \nu_{2a} + 4 \nu_{2b} + 4 \nu_{4a} + \alpha) \geq 0; \\ 
\mu_{2}(u,\chi_{2},*) & = \textstyle \frac{1}{4} (4 \nu_{2a} - 4 \nu_{2b} - 4 \nu_{4a} + \alpha) \geq 0; \\ 
\mu_{0}(u,\chi_{4},*) & = \textstyle \frac{1}{4} (10 \nu_{2a} - 6 \nu_{2b} + 2 \nu_{4a} + \beta) \geq 0; \\ 
\mu_{2}(u,\chi_{4},*) & = \textstyle \frac{1}{4} (-10 \nu_{2a} + 6 \nu_{2b} - 2 \nu_{4a} + \beta) \geq 0, \\ 
\end{split}
\]
\[
\begin{split}
\text{where} \quad (\alpha,\beta)=
\tiny{
\begin{cases}
(12,26) \quad \text{when} \quad \chi(u^2)=\chi(2a); \\
(16,18) \quad \text{when} \quad \chi(u^2)=\chi(2b); \\
(24,2) \quad \text{when} \quad \chi(u^2)=-2\chi(2a)+3\chi(2b); \\
(8,34) \quad \text{when} \quad \chi(u^2)= 2\chi(2a)-\chi(2b); \\
(4,42) \quad \text{when} \quad \chi(u^2)= 3\chi(2a)-2\chi(2b); \\
(20,10) \quad \text{when} \quad \chi(u^2)= -\chi(2a)+2\chi(2b). \\
\end{cases}}
\end{split}
\]
Additionally, we need to consider the following case-dependent inequalities:
\[
\begin{split}
\mu_{0}(u,\chi_{2},3) & = \textstyle \frac{1}{4} (-6 \nu_{2a} + 2 \nu_{2b} + 2 \nu_{4a} + 10) \geq 0 \quad
\text{for} \; \chi(u^2)=\chi(2a); \\ 
\mu_{0}(u,\chi_{8},*) & = \textstyle \frac{1}{4} (-20 \nu_{2a} - 4 \nu_{2b} + 4 \nu_{4a} + 68) \geq 0 \quad
\text{for} \; \chi(u^2)=\chi(2b); \\ 
\mu_{2}(u,\chi_{7},*) & = \textstyle \frac{1}{4} (-30 \nu_{2a} + 2 \nu_{2b} - 6 \nu_{4a} + 30) \geq 0 \quad
\text{for} \; \chi(u^2)=-2\chi(2a)+3\chi(2b); \\ 
\mu_{0}(u,\chi_{4},5) & = \textstyle \frac{1}{4} (18 \nu_{2a} + 2 \nu_{2b} + 2 \nu_{4a} + 26) \geq 0 \quad
\text{for} \; \chi(u^2)=-2\chi(2a)+3\chi(2b); \\ 
\mu_{0}(u,\chi_{8},*) & = \textstyle \frac{1}{4} (-20 \nu_{2a} - 4 \nu_{2b} + 4 \nu_{4a} + 52) \geq 0 \quad
\text{for} \; \chi(u^2)= 2\chi(2a)-\chi(2b); \\ 
\mu_{0}(u,\chi_{2},3) & = \textstyle \frac{1}{4} (-6 \nu_{2a} + 2 \nu_{2b} + 2 \nu_{4a} + 2) \geq 0 \quad
\text{for} \; \chi(u^2)= 3\chi(2a)-2\chi(2b); \\ 
\mu_{0}(u,\chi_{8},*) & = \textstyle \frac{1}{4} (-20 \nu_{2a} - 4 \nu_{2b} + 4 \nu_{4a} + 76) \geq 0 \quad
\text{for} \; \chi(u^2)= -\chi(2a)+2\chi(2b); \\ 
\mu_{0}(u,\chi_{10},*) & = \textstyle \frac{1}{4} (20 \nu_{2a} + 12 \nu_{2b} - 4 \nu_{4a} + 92) \geq 0 \quad
\text{for} \; \chi(u^2)= -\chi(2a)+2\chi(2b). \\ 
\end{split}
\]
Solving these systems and applying Proposition \ref{P:6} to the obtained solutions, we get only
fifteen integer solutions $(\nu_{2a},\nu_{2b},\nu_{4a})$ as listed in part (v) of Theorem \ref{T:2}.

\noindent $\bullet$ Let $|u|=5$.
By (\ref{E:1}) and Proposition \ref{P:4} we get
$\nu_{5a}+\nu_{5b}+\nu_{5c}+\nu_{5d}=1$. Put $t=20 \nu_{5a} + 20 \nu_{5b}$.
By Proposition \ref{P:1} we get the system of inequalities:
\[
\begin{split}
\mu_{0}(u,\chi_{10},*) & = \textstyle \frac{1}{5} ( t + 90) \geq 0; \quad 
\mu_{0}(u,\chi_{12},*)   = \textstyle \frac{1}{5} (-t + 160) \geq 0; \\ 
& \mu_{1}(u,\chi_{8},*)   = \textstyle \frac{1}{5} (-15 \nu_{5a} + 10 \nu_{5b} + 70) \geq 0; \\ 
& \mu_{2}(u,\chi_{8},*)   = \textstyle \frac{1}{5} (10 \nu_{5a} - 15 \nu_{5b} + 70) \geq 0, \\ 
\end{split}
\]
from which we can derive 71 possible pairs $(\nu_{5a},\nu_{5b})$. From the inequalities:
\[
\begin{split}
\mu_{0}(u,\chi_{7},2) & = \textstyle \frac{1}{5} (16 \nu_{5a} + 16 \nu_{5b} - 4 \nu_{5c} - 4 \nu_{5d} + 64) \geq 0; \\ 
\mu_{0}(u,\chi_{6},2) & = \textstyle \frac{1}{5} (-16 \nu_{5a} - 16 \nu_{5b} + 4 \nu_{5c} + 4 \nu_{5d} + 36) \geq 0, \\ 
\end{split}
\]
if follows that $t = 4 \nu_{5a} + 4 \nu_{5b} -  \nu_{5c} -  \nu_{5d} \in \{ -16, -11, -6, -1, 4, 9 \}$.
Taking into account that $\nu_{5a}+\nu_{5b}+\nu_{5c}+\nu_{5d}=1$
and considering the additional inequality
$$
\mu_{0}(u,\chi_{2},2) = \textstyle \frac{1}{5} (4 \nu_{5a} + 4 \nu_{5b} - 6 \nu_{5c} - 6 \nu_{5d} + 6) \geq 0, 
$$
it is easy to check that there remain only 16 possibilities for
$(\nu_{5a},\nu_{5b},\nu_{5c}+\nu_{5d})$:
\[
\begin{split}
\{ &
( -2, 2, 1 ), ( -2, 3, 0 ), ( -1, 1, 1 ), ( -1, 2, 0 ), ( -1, 3, -1 ), ( 0, 0, 1 ), ( 0, 1, 0 ), ( 0, 2, -1 ), \\
& ( 1, -1, 1 ), ( 1, 0, 0 ), ( 1, 1, -1 ), ( 2, -2, 1 ), ( 2, -1, 0 ), ( 2, 0, -1 ), ( 3, -2, 0 ), ( 3, -1, -1 ) \}. \\
\end{split}
\]
Finally, using the inequalities
\[
\begin{split}
\mu_{1}(u,\chi_{2},2) & = \textstyle \frac{1}{5} (-6 \nu_{5a} + 4 \nu_{5b} -  \nu_{5c} + 4 \nu_{5d} + 6) \geq 0; \\
\mu_{2}(u,\chi_{2},2) & = \textstyle \frac{1}{5} (4 \nu_{5a} - 6 \nu_{5b} + 4 \nu_{5c} -  \nu_{5d} + 6) \geq 0; \\
\mu_{1}(u,\chi_{4},2) & = \textstyle \frac{1}{5} (-9 \nu_{5a} + 6 \nu_{5b} +  \nu_{5c} - 4 \nu_{5d} + 14) \geq 0; \\
\mu_{2}(u,\chi_{4},2) & = \textstyle \frac{1}{5} (6 \nu_{5a} - 9 \nu_{5b} - 4 \nu_{5c} +  \nu_{5d} + 14) \geq 0; \\
\mu_{1}(u,\chi_{7},2) & = \textstyle \frac{1}{5} (6 \nu_{5a} - 14 \nu_{5b} + 6 \nu_{5c} - 4 \nu_{5d} + 64) \geq 0; \\ 
\mu_{2}(u,\chi_{7},2) & = \textstyle \frac{1}{5} (-14 \nu_{5a} + 6 \nu_{5b} - 4 \nu_{5c} + 6 \nu_{5d} + 64) \geq 0, \\ 
\end{split}
\]
we obtain only ten integer solutions listed in part (vi) of Theorem \ref{T:2}.

\noindent $\bullet$ Let $|u|=7$. Since
there is only one conjugacy class in $G$
consisting of elements or order 7, this case
follows immediately from Proposition \ref{P:4}.

\noindent $\bullet$ Let $|u|=8$.
By (\ref{E:1}) and Proposition \ref{P:4} we get
$\nu_{2a}+\nu_{2b}+\nu_{4a}+\nu_{8a}=1$.
Because $|u^2|=4$ and $|u^4|=2$,
we need to consider 90 cases defined  by parts (iii) and (v) of Theorem \ref{T:2}.
First, in 45 of these cases, given in the following table, we have no units of order 8
because $\mu_1(u,\chi_2,*)$ is not an integer:
$$
\tiny{
\begin{array}{|c|c|c|c|}\hline
                              &                         &                         & \\
\chi(u^2)                     & \chi(u^4)= \chi(2b)     &
                                \chi(u^4)= -2\chi(2a)+3\chi(2b)  &
                                \chi(u^4)= 2\chi(2a)-\chi(2b) \\
                              &                         &                         & \\ \hline
\chi(4a)                      &                         &                         &                         \\
-2\chi(2a)-2\chi(2b)+5\chi(4a)&                         &                         &                         \\
-\chi(2a)-3\chi(2b)+5\chi(4a) &                         &                         &                         \\
-\chi(2a)-\chi(2b)+3\chi(4a)  &                         &                         &                         \\
-\chi(2a)+\chi(2b)+\chi(4a)   &                         &                         &                         \\
-4\chi(2b)+5\chi(4a)          &                         &                         &                         \\
-2\chi(2b)+3\chi(4a)          &                         &                         &                         \\
2\chi(2b)-\chi(4a)            & \mu_1(u,\chi_2,*)= \textstyle\frac{3}{2} &
                                \mu_1(u,\chi_2,*)= \textstyle\frac{1}{2} &
                                \mu_1(u,\chi_2,*)= \textstyle\frac{5}{2} \\
4\chi(2b)-3\chi(4a)           &                         &                         &                         \\
\chi(2a)-3\chi(2b)+3\chi(4a)  &                         &                         &                         \\
\chi(2a)-\chi(2b)+\chi(4a)    &                         &                         &                         \\
\chi(2a)+\chi(2b)-\chi(4a)    &                         &                         &                         \\
\chi(2a)+3\chi(2b)-3\chi(4a)  &                         &                         &                         \\
2\chi(2a)-\chi(4a)            &                         &                         &                         \\
2\chi(2a)+2\chi(2b)-3\chi(4a) &                         &                         &                         \\ \hline
\end{array}}
$$
Now put $t_1=8 \nu_{2a} - 8 \nu_{2b} - 8 \nu_{4a}$ and
$t_2=20 \nu_{2a} - 12 \nu_{2b} + 4 \nu_{4a} - 4 \nu_{8a}$.
Then, when $\chi(u^{4}) = 3\chi(2a)-2\chi(2b)$ and $\chi(u^{2})$ is equal to $2\chi(2a)-\chi(4a)$
or $2\chi(2a)+2\chi(2b)-3\chi(4a)$, we obtain the system of inequalities:
$$
\mu_{0}(u,\chi_{2},*) = \textstyle \frac{1}{8} (-t_1 - 8) \geq 0; \quad 
\mu_{4}(u,\chi_{2},*) = \textstyle \frac{1}{8} (t_1 - 8) \geq 0, 
$$
which have no integer solutions.
Also, there is no solution for the system
$$
\mu_{0}(u,\chi_{4},*) = \textstyle \frac{1}{8} (t_2 + \alpha) \geq 0; \quad 
\mu_{4}(u,\chi_{4},*) = \textstyle \frac{1}{8} (-t_2 + \alpha) \geq 0, 
$$
where $\alpha=-4$ for $(\chi(u^{4}), \chi(u^{2}))$ in the set
\[
\begin{split}
\{ \; & (\chi(2a), \; 4\chi(2b)-3\chi(4a)), \;
(-\chi(2a)+2\chi(2b), \; -\chi(2a)+\chi(2b)+\chi(4a)), \\ &
(-\chi(2a)+2\chi(2b), \; 2\chi(2b)-\chi(4a)), \;
(-\chi(2a)+2\chi(2b), \; \chi(2a)+3\chi(2b)-3\chi(4a)) \; \}
\end{split}
\]
and $\alpha=-20$ \quad for $(\chi(u^{4}),\chi(u^{2}))=(-\chi(2a)+2\chi(2b),4\chi(2b)-3\chi(4a))$.

In the remaining 38 cases we first consider the following system of inequalities:
\[
\begin{split}
\mu_{0}(u,\chi_{2},*) & = \textstyle \frac{1}{8} (-t_1 + \alpha_1) \geq 0; \quad 
\mu_{4}(u,\chi_{2},*)  = \textstyle \frac{1}{8} (t_1 + \alpha_1) \geq 0; \\ 
\mu_{0}(u,\chi_{4},*) & = \textstyle \frac{1}{8} (t_2 + \alpha_2) \geq 0; \quad 
\mu_{4}(u,\chi_{4},*)  = \textstyle \frac{1}{8} (-t_2 + \alpha_2) \geq 0; \\ 
\mu_{0}(u,\chi_{7},*) & = \textstyle \frac{1}{8} (t_3 + \alpha_3) \geq 0; \quad 
\mu_{4}(u,\chi_{7},*)  = \textstyle \frac{1}{8} (-t_3 + \alpha_3) \geq 0, \\ 
\end{split}
\]
where $t_3=60 \nu_{2a} - 4 \nu_{2b} + 12 \nu_{4a} + 4 \nu_{8a}$ and
the tuples $(\alpha_1,\alpha_2,\alpha_3)$ are given below:
$$
\tiny{
\begin{array}{|c|c|c|c|}\hline
& \chi(u^{4})
& \chi(u^{2})
& (\alpha_1,\alpha_2,\alpha_3) \\ \hline
1
&
& \chi(4a)
& ( 16, 28, 84 ) \\
2
&
& -2\chi(2a)-2\chi(2b)+5\chi(4a)
& ( 32, 28, 52 ) \\
3
&
& -\chi(2a)-3\chi(2b)+5\chi(4a)
& ( 24, 44, 84 ) \\
4
&
& -\chi(2a)-\chi(2b)+3\chi(4a)
& ( 24, 28, 68 ) \\
5
&
& -\chi(2a)+\chi(2b)+\chi(4a)
& ( 24, 12, 52 ) \\
6
&
& -4\chi(2b)+5\chi(4a)
& ( 16, 60, 116 ) \\
7
& \chi(2a)
& -2\chi(2b)+3\chi(4a)
& ( 16, 44, 100 ) \\
8
&
& 2\chi(2b)-\chi(4a)
& ( 16, 12, 68 ) \\
9
&
& \chi(2a)-3\chi(2b)+3\chi(4a)
& ( 8, 60, 132 ) \\
10
&
& \chi(2a)-\chi(2b)+\chi(4a)
& ( 8, 44, 116 ) \\
11
&
& \chi(2a)+\chi(2b)-\chi(4a)
& ( 8, 28, 100 ) \\
12
&
& \chi(2a)+3\chi(2b)-3\chi(4a)
& ( 8, 12, 84 ) \\
13
&
& 2\chi(2a)-\chi(4a)
& (0, 44, 132) \\
14
&
& 2\chi(2a)+2\chi(2b)-3\chi(4a)
& ( 0, 28, 116 ) \\ \hline
15
&
& \chi(4a)
& ( 8, 44, 116 ) \\
16
&
& -2\chi(2a)-2\chi(2b)+5\chi(4a)
& ( 24, 44, 84 ) \\
17
&
& -\chi(2a)-3\chi(2b)+5\chi(4a)
& ( 16, 60, 116 ) \\
18
&
& -\chi(2a)-\chi(2b)+3\chi(4a)
& ( 16, 44, 100 ) \\
19
&
& -\chi(2a)+\chi(2b)+\chi(4a)
& ( 16, 28, 84 ) \\
20
&
& -4\chi(2b)+5\chi(4a)
& ( 8, 76, 148 ) \\
21
& 3\chi(2a)-2\chi(2b)
& -2\chi(2b)+3\chi(4a)
& ( 8, 60, 132 ) \\
22
&
& 2\chi(2b)-\chi(4a)
& ( 8, 28, 100 ) \\
23
&
& 4\chi(2b)-3\chi(4a)
& ( 8, 12, 84 ) \\
24
&
& \chi(2a)-3\chi(2b)+3\chi(4a)
& ( 0, 76, 164 ) \\
25
&
& \chi(2a)-\chi(2b)+\chi(4a)
& ( 0, 60, 148 ) \\
26
&
& \chi(2a)+\chi(2b)-\chi(4a)
& ( 0, 44, 132 ) \\
27
&
& \chi(2a)+3\chi(2b)-3\chi(4a)
& ( 0, 28, 116 ) \\ \hline
28
&
& \chi(4a)
& ( 24, 12, 52 ) \\
29
&
& -2\chi(2a)-2\chi(2b)+5\chi(4a)
& ( 40, 12, 20 ) \\
30
&
& -\chi(2a)-3\chi(2b)+5\chi(4a)
& ( 32, 28, 52 ) \\
31
&
& -\chi(2a)-\chi(2b)+3\chi(4a)
& ( 32, 12, 36 ) \\
32
&
& -4\chi(2b)+5\chi(4a)
& ( 24, 44, 84 ) \\
33
& -\chi(2a)+2\chi(2b)
& -2\chi(2b)+3\chi(4a)
& ( 24, 28, 68 ) \\
34
&
& \chi(2a)-3\chi(2b)+3\chi(4a)
& ( 16, 44, 100 ) \\
35
&
& \chi(2a)-\chi(2b)+\chi(4a)
& ( 16, 28, 84 ) \\
36
&
& \chi(2a)+\chi(2b)-\chi(4a)
& ( 16, 12, 68 ) \\
37
&
& 2\chi(2a)-\chi(4a)
& ( 8, 28, 100 ) \\
38
&
& 2\chi(2a)+2\chi(2b)-3\chi(4a)
& ( 8, 12, 84 ) \\
\hline
\end{array}}
$$
In all 38 cases these inequalities allow us to compute admissible solutions,
using the technique explained in detail in the proof of Theorem \ref{T:1} in the
case of units of order 10. Having done this, we need to consider additional inequalities
and apply Proposition \ref{P:6} to reduce the number of solutions or, possibly,
eliminate all of them.

In cases 1, 19 and 35 we use the system of inequalities:
\[
\begin{split}
\mu_{0}(u,\chi_{2},3) & = \textstyle \frac{1}{8} (-12 \nu_{2a} + 4 \nu_{2b} + 4 \nu_{4a} - 4 \nu_{8a} + 12) \geq 0; \\ 
\mu_{0}(u,\chi_{7},3) & = \textstyle \frac{1}{8} (-28 \nu_{2a} - 12 \nu_{2b} + 4 \nu_{4a} + 4 \nu_{8a} + 52) \geq 0; \\ 
\mu_{4}(u,\chi_{7},3) & = \textstyle \frac{1}{8} (28 \nu_{2a} + 12 \nu_{2b} - 4 \nu_{4a} - 4 \nu_{8a} + 52) \geq 0, \\ 
\end{split}
\]
to obtain in all three cases the same set of solutions $\{$
$( 0, 2, 0, -1 )$, $( 0, -2, 2, 1 )$, $( 0, 0, 2, -1 )$, $( 0, 0, 0, 1 )$,
$( -1, -1, 2, 1 )$, $( -1, -1, 0, 3 )$, $( 1, 1, 2, -3 )$, $( 1, 1, 0, -1 )$ $\}$.

In cases 2, 3, 4, 5, 7, 16, 18, 28, 30, 31, 32, 33 and 34 we use the system:
\[
\begin{split}
\mu_{4}(u,\chi_{10},*) & = \textstyle \frac{1}{8} (-40 \nu_{2a} - 24 \nu_{2b} + 8 \nu_{4a} + \beta_1) \geq 0; \\ 
\mu_{0}(u,\chi_{2},3) & = \textstyle \frac{1}{8} (-12 \nu_{2a} + 4 \nu_{2b} + 4 \nu_{4a} - 4 \nu_{8a} + \beta_2) \geq 0; \\ 
\mu_{0}(u,\chi_{10},7) & = \textstyle \frac{1}{8} (36 \nu_{2a} + 20 \nu_{2b} - 12 \nu_{4a} - 4 \nu_{8a} + \beta_3) \geq 0, \\ 
\end{split}
\]
and the following table describes tuples $(\beta_1,\beta_2,\beta_3)$ and solutions for each case:
$$
\small{
\begin{array}{|l|c|c|}\hline
\text{Cases} & (\beta_1,\beta_2,\beta_3) & (\nu_{2a}, \; \nu_{2b}, \; \nu_{4a}, \; \nu_{8a}) \\ \hline
2,30   &(16,28,12)& ( 0, 0, -2, 3 ), \; ( 0, 0, 0, 1 ) \\ \hline
3,16,32&(24,20,20)& ( 0, 0, 2, -1 ), \; ( 0, 0, -2, 3 ), \; ( 0, 0, 0, 1 ), \\
       &          & ( 1, -1, 2, -1 ), \; ( 1, -1, 0, 1 ), \; ( -1, 1, 0, 1 ) \\ \hline
4,33   &(56,20,52)& ( 0, 2, 0, -1 ), \; ( 0, -2, 0, 3 ), \; ( 0, 0, 2, -1 ), \; ( 0, 0, -2, 3 ), \\
       &          & ( 0, 0, 0, 1 ), \; ( 1, -1, 0, 1 ), \; ( 1, 1, 2, -3 ) \\ \hline
5,28   &(88,20,84)& ( 0, -2, 0, 3 ), \; ( 0, 0, 2, -1 ), \; ( 0, 0, 0, 1 ), \\
       &          & ( -1, -1, 2, 1 ), \; ( 1, 1, -2, 1 ), \; ( 1, 1, 0, -1 ) \\ \hline
7,18,34&(64,12,60)& ( 0, 2, -2, 1 ), \; ( 0, 2, 0, -1 ), \; ( 0, 0, 2, -1 ), \; ( 0, 0, 0, 1 ), \\
       &          & ( 1, -1, 2, -1 ), \; ( 1, 1, 2, -3 ), \; ( 1, 1, 0, -1 ), \; ( -1, 1, 0, 1 ) \\ \hline
31     &(48,28,44)& ( 0, 2, 2, -3 ), \; ( 0, 0, 2, -1 ), \; ( 0, 0, 0, 1 ) \\ \hline
\end{array}}
$$
In cases 6 and 17 we use the system of inequalities
\[
\begin{split}
\mu_{4}(u,\chi_{10},*) & = \textstyle \frac{1}{8} (-40 \nu_{2a} - 24 \nu_{2b} + 8 \nu_{4a} + 32) \geq 0; \\ 
\mu_{0}(u,\chi_{2},3) & = \textstyle \frac{1}{8} (-12 \nu_{2a} + 4 \nu_{2b} + 4 \nu_{4a} - 4 \nu_{8a} + 12) \geq 0; \\ 
\mu_{0}(u,\chi_{10},7) & = \textstyle \frac{1}{8} (36 \nu_{2a} + 20 \nu_{2b} - 12 \nu_{4a} - 4 \nu_{8a} + 28) \geq 0; \\ 
\mu_{0}(u,\chi_{12},7) & = \textstyle \frac{1}{8} (-16 \nu_{2a} + 16 \nu_{2b} - 16 \nu_{4a} + 48) \geq 0; \\ 
\mu_{4}(u,\chi_{12},7) & = \textstyle \frac{1}{8} (16 \nu_{2a} - 16 \nu_{2b} + 16 \nu_{4a} + 48) \geq 0, \\ 
\end{split}
\]
to obtain in both cases the same set of solutions
$\{$ $(0, 0, 2, -1)$, $(0, 0, 0, 1)$, $(-1, 1, 0, 1) \}$.

In cases 8 and 36 we use the system of inequalities:
\[
\begin{split}
\mu_{4}(u,\chi_{8},*) & = \textstyle \frac{1}{8} (40 \nu_{2a} + 8 \nu_{2b} - 8 \nu_{4a} + 48) \geq 0; \\ 
\mu_{0}(u,\chi_{2},3) & = \textstyle \frac{1}{8} (-12 \nu_{2a} + 4 \nu_{2b} + 4 \nu_{4a} - 4 \nu_{8a} + 12) \geq 0; \\ 
\mu_{0}(u,\chi_{7},3) & = \textstyle \frac{1}{8} (-28 \nu_{2a} - 12 \nu_{2b} + 4 \nu_{4a} + 4 \nu_{8a} + 36) \geq 0, \\ 
\end{split}
\]
to obtain in both cases the same set of solutions
$\{ (0, 0, 2, -1), (0, 0, 0, 1) \}$.

In cases 9, 20 and 21 we use the system
\[
\begin{split}
\mu_{4}(u,\chi_{10},*) & = \textstyle \frac{1}{8} (-40 \nu_{2a} - 24 \nu_{2b} + 8 \nu_{4a} + \beta_1) \geq 0; \\ 
\mu_{0}(u,\chi_{2},3) & = \textstyle \frac{1}{8} (-12 \nu_{2a} + 4 \nu_{2b} + 4 \nu_{4a} - 4 \nu_{8a} + 4) \geq 0; \\ 
\mu_{0}(u,\chi_{10},7) & = \textstyle \frac{1}{8} (36 \nu_{2a} + 20 \nu_{2b} - 12 \nu_{4a} - 4 \nu_{8a} + \beta_2) \geq 0; \\ 
\mu_{0}(u,\chi_{12},7) & = \textstyle \frac{1}{8} (-16 \nu_{2a} + 16 \nu_{2b} - 16 \nu_{4a} + \beta_3) \geq 0; \\ 
\mu_{4}(u,\chi_{12},7) & = \textstyle \frac{1}{8} (16 \nu_{2a} - 16 \nu_{2b} + 16 \nu_{4a} + \beta_3) \geq 0, \\ 
\end{split}
\]
where $(\beta_1,\beta_2,\beta_3)$ is equal to $(72,68,64)$ in cases 9 and 21,
and to $(40,36,32)$ in case 20. This leads to the solutions
$( 0, 2, -2, 1 )$, \; $( 0, -2, 2, 1 )$, \; $( 2, 0, 2, -3 )$, \; $( 0, 0, 0, 1 )$, \;
$( 1, -1, 2, -1 )$, \; $( -1, -1, 0, 3 )$, \; $( -1, 1, -2, 3 )$, \; $( 1, 1, 0, -1)$
in cases 9 and 21 and to the unique solution $(0,0,0,1)$ in case 20.

In cases 10, 11, 12, 15, 22, 23, 37 and 38 we use the system
\[
\begin{split}
\mu_{0}(u,\chi_{8},*) & = \textstyle \frac{1}{8} (-40 \nu_{2a} - 8 \nu_{2b} + 8 \nu_{4a} + \beta_1) \geq 0; \\ 
\mu_{0}(u,\chi_{2},3) & = \textstyle \frac{1}{8} (-12 \nu_{2a} + 4 \nu_{2b} + 4 \nu_{4a} - 4 \nu_{8a} + 4) \geq 0; \\ 
\mu_{4}(u,\chi_{7},3) & = \textstyle \frac{1}{8} (28 \nu_{2a} + 12 \nu_{2b} - 4 \nu_{4a} - 4 \nu_{8a} + \beta_2) \geq 0, \\ 
\end{split}
\]
and the following table describes tuples $(\beta_1,\beta_2)$ and solutions for each case:
$$
\small{
\begin{array}{|l|c|c|}\hline
\text{Cases} & (\beta_1,\beta_2) & (\nu_{2a}, \; \nu_{2b}, \; \nu_{4a}, \; \nu_{8a}) \\ \hline
10,15    &(48,44)& (0, 2, -2, 1), \; (0, -2, 2, 1), \; (0, 0, 0, 1), \; (1, -1, 2, -1), \; (1, 1, 0, -1) \\
11,22,37 &(32,28)& (0,0,0,1) \\
12,23,38 &(16,12)& (0,0,0,1) \\ \hline
\end{array}}
$$
Finally, in cases 13, 14, 24, 25, 26, 27 we use the additional inequality
$$
\mu_{0}(u,\chi_{2},3) = \textstyle \frac{1}{8} (-12 \nu_{2a} + 4 \nu_{2b} + 4 \nu_{4a} - 4 \nu_{8a} - 4) \geq 0, 
$$
and in the case 29 we use two additional inequalities:
\[
\begin{split}
\mu_{4}(u,\chi_{10},*) & = \textstyle \frac{1}{8} (-40 \nu_{2a} - 24 \nu_{2b} + 8 \nu_{4a} + 8) \geq 0; \\ 
\mu_{0}(u,\chi_{10},7) & = \textstyle \frac{1}{8} (36 \nu_{2a} + 20 \nu_{2b} - 12 \nu_{4a} - 4 \nu_{8a} + 4) \geq 0, \\ 
\end{split}
\]
to show that in these cases we have no solutions.

Now the union of the solutions obtained above gives us part (vii) of Theorem \ref{T:2}.

\noindent $\bullet$ Let $|u|= 15$. By
(\ref{E:1}) and Proposition \ref{P:4} we get $$
\nu_{3a}+\nu_{3b}+\nu_{5a}+\nu_{5b}+\nu_{5c}+\nu_{5d}+\nu_{15a}+\nu_{15b}=1.
$$ We need to consider 30 cases defined by parts
(iv) and (vi) of Theorem \ref{T:2}. Only in two
cases we will get a system of inequalities that
has solutions. Furthermore, each time these
solutions are trivial, so they will give the proof for order $15$.

First, in eight cases we obtain the following system of inequalities
that has no solutions such that all $\mu_{i}(u,\chi_{j},*)$ are non-negative integers:
\[
\begin{split}
\mu_{0}(u,\chi_{6},*) & = \textstyle \frac{1}{15} (72 \nu_{3a} - 32 \nu_{5a} - 32 \nu_{5b} + 8 \nu_{5c} + 8 \nu_{5d} - 8 \nu_{15a} - 8 \nu_{15b} + \alpha) \geq 0; \\
\mu_{5}(u,\chi_{6},*) & = \textstyle \frac{1}{15} (-36 \nu_{3a} + 16 \nu_{5a} + 16 \nu_{5b} - 4 \nu_{5c} - 4 \nu_{5d} + 4 \nu_{15a} + 4 \nu_{15b} + \beta) \geq 0, \\
\end{split}
\]
where the tuples $(\alpha,\beta)$ are given in the following table:
$$
\small{
\begin{array}{|c|c|c|}\hline
\chi(u^5)           & \chi(u^3)                  & (\alpha,\beta) \\ \hline
\chi(3a)            & \chi(5a)+\chi(5b)-\chi(5d) & (18,-9)        \\        
                    & \chi(5a)+\chi(5b)-\chi(5c) &                \\ \hline 
                    & \chi(5a) , \quad \chi(5b)                   & (2,29)         \\        
-\chi(3a)+2\chi(3b) & \chi(5a)-\chi(5c)+\chi(5d) &                \\        
                    & \chi(5b)+\chi(5c)-\chi(5d) &                \\ \hline 
-\chi(3a)+2\chi(3b) & \chi(5a)+\chi(5b)-\chi(5d) & (-18,9)        \\        
                    & \chi(5b)+\chi(5c)-\chi(5d) &                \\ \hline 
\end{array}}
$$
In four other cases we obtain the following system of inequalities
that has no solution such that all $\mu_{i}(u,\chi_{j},*)$ are non-negative integers:
\[
\small{
\begin{split}
\mu_{0}(u,\chi_{2},*) &=\textstyle\frac{1}{15}(40 \nu_{3a} - 8 \nu_{3b} + 12 \nu_{5a} + 12 \nu_{5b} + 12 \nu_{5c} + 12 \nu_{5d} + 6) \geq 0; \\ 
\mu_{1}(u,\chi_{2},*) &=\textstyle\frac{1}{15}(5 \nu_{3a} -  \nu_{3b} - 6 \nu_{5a} + 9 \nu_{5b} + 4 \nu_{5c} -  \nu_{5d} + \alpha_1) \geq 0; \\ 
\mu_{3}(u,\chi_{2},*) &=\textstyle\frac{1}{15} (-10 \nu_{3a} + 2 \nu_{3b} - 18 \nu_{5a} + 12 \nu_{5b} + 2 \nu_{5c} - 8 \nu_{5d} + \alpha_7) \geq 0; \\ 
\mu_{5}(u,\chi_{2},*) &=\textstyle\frac{1}{15}(-20 \nu_{3a} + 4 \nu_{3b} - 6 \nu_{5a} - 6 \nu_{5b} - 6 \nu_{5c} - 6 \nu_{5d} + 27) \geq 0; \\ 
\mu_{6}(u,\chi_{2},*) &=\textstyle\frac{1}{15}(-10 \nu_{3a} + 2 \nu_{3b} + 12 \nu_{5a} - 18 \nu_{5b} - 8 \nu_{5c} + 2 \nu_{5d} + \alpha_2) \geq 0; \\ 
\mu_{0}(u,\chi_{4},*) &=\textstyle\frac{1}{15}(24 \nu_{3a} + 28 \nu_{5a} + 28 \nu_{5b} + 8 \nu_{5c} + 8 \nu_{5d} + 4 \nu_{15a} + 4 \nu_{15b} + 19) \geq 0; \\ 
\mu_{1}(u,\chi_{4},*) &=\textstyle\frac{1}{15}(3 \nu_{3a} + 6 \nu_{5a} +  \nu_{5b} - 4 \nu_{5c} + 6 \nu_{5d} + 3 \nu_{15a} - 2 \nu_{15b} + \alpha_3) \geq 0; \\ 
\mu_{5}(u,\chi_{4},*) &=\textstyle\frac{1}{15}(-12 \nu_{3a} - 14 \nu_{5a} - 14 \nu_{5b} - 4 \nu_{5c} - 4 \nu_{5d} - 2 \nu_{15a} - 2 \nu_{15b} + 28) \geq 0; \\ 
\mu_{6}(u,\chi_{4},*) &=\textstyle\frac{1}{15}(-6 \nu_{3a} - 12 \nu_{5a} - 2 \nu_{5b} + 8 \nu_{5c} - 12 \nu_{5d} - 6 \nu_{15a} + 4 \nu_{15b} + \alpha_4) \geq 0; \\ 
\mu_{0}(u,\chi_{6},*) &=\textstyle\frac{1}{15}(72 \nu_{3a} - 32 \nu_{5a} - 32 \nu_{5b} + 8 \nu_{5c} + 8 \nu_{5d} - 8 \nu_{15a} - 8 \nu_{15b} + 22) \geq 0; \\ 
\mu_{3}(u,\chi_{6},*) &=\textstyle\frac{1}{15}(-18 \nu_{3a} + 8 \nu_{5a} + 8 \nu_{5b} - 2 \nu_{5c} - 2 \nu_{5d} + 2 \nu_{15a} + 2 \nu_{15b} + 17) \geq 0; \\ 
\mu_{0}(u,\chi_{11}, *) & =\textstyle\frac{1}{15}(-72 \nu_{3a} + 8 \nu_{5a} + 8 \nu_{5b} + 8 \nu_{5c} + 8 \nu_{5d} + 8 \nu_{15a} + 8 \nu_{15b} + 148) \geq 0; \\ 
\mu_{5}(u,\chi_{11}, *) & =\textstyle\frac{1}{15}(36 \nu_{3a} - 4 \nu_{5a} - 4 \nu_{5b} - 4 \nu_{5c} - 4 \nu_{5d} - 4 \nu_{15a} - 4 \nu_{15b} + 121) \geq 0; \\ 
\mu_{1}(u,\chi_{16}, *) & =\textstyle\frac{1}{15}(8 \nu_{3a} -  \nu_{3b} + 9 \nu_{5a} - 11 \nu_{5b} + 4 \nu_{5c} - 6 \nu_{5d} + 3 \nu_{15a} - 2 \nu_{15b} + \alpha_5) \geq 0;\\
\mu_{6}(u,\chi_{16}, *) & =\textstyle\frac{1}{15}(-16 \nu_{3a} + 2 \nu_{3b} - 18 \nu_{5a} + 22 \nu_{5b} - 8 \nu_{5c} \\
& \hspace{150pt} + 12 \nu_{5d} - 6 \nu_{15a} + 4 \nu_{15b} + \alpha_6)\geq 0,\\
\end{split}}
\]
where the tuples $(\alpha_1,\alpha_2,\alpha_3,\alpha_4,\alpha_5,\alpha_6,\alpha_7)$ are given in the following table:
$$
\small{
\begin{array}{|c|c|c|}\hline
\qquad \quad \chi(u^5) \quad \qquad & \chi(u^3)           & (\alpha_1,\alpha_2,\alpha_3,\alpha_4,\alpha_5,\alpha_6,\alpha_7)  \\ \hline
-\chi(3a)+2\chi(3b) & \chi(5c)            & (22,1,18,9,240,210,-4)  \\ \hline  
-\chi(3a)+2\chi(3b) & \chi(5d)            & (17,-4,28,19,230,200,1) \\ \hline  
-\chi(3a)+2\chi(3b) & 2\chi(5c)-\chi(5d)  & (27,6,8,-1,250,220,-9)  \\ \hline  
-\chi(3a)+2\chi(3b) & -\chi(5c)+2\chi(5d) & (12,-9,38,29,220,190,6) \\ \hline  
\end{array}}
$$
In the remaining 18 cases we first consider the following system of inequalities:
\[
\small{
\begin{split}
\mu_{0}(u,\chi_{2}, *) & =\textstyle\frac{1}{15}(40 \nu_{3a} - 8 \nu_{3b} + 12 \nu_{5a} + 12 \nu_{5b} + 12 \nu_{5c} + 12 \nu_{5d} + \alpha_1)\geq 0;\\
\mu_{1}(u,\chi_{2}, *) & =\textstyle\frac{1}{15}(5 \nu_{3a} -  \nu_{3b} - 6 \nu_{5a} + 9 \nu_{5b} + 4 \nu_{5c} -  \nu_{5d} + \alpha_2)\geq 0;\\
\mu_{5}(u,\chi_{2}, *) & =\textstyle\frac{1}{15}(-20 \nu_{3a} + 4 \nu_{3b} - 6 \nu_{5a} - 6 \nu_{5b} - 6 \nu_{5c} - 6 \nu_{5d} + \alpha_3)\geq 0;\\
\mu_{6}(u,\chi_{2}, *) & =\textstyle\frac{1}{15}(-10 \nu_{3a} + 2 \nu_{3b} + 12 \nu_{5a} - 18 \nu_{5b} - 8 \nu_{5c} + 2 \nu_{5d} + \alpha_4)\geq 0;\\
\mu_{0}(u,\chi_{4}, *) & =\textstyle\frac{1}{15}(24 \nu_{3a} + 28 \nu_{5a} + 28 \nu_{5b} + 8 \nu_{5c} + 8 \nu_{5d} + 4 \nu_{15a} + 4 \nu_{15b} + \alpha_5)\geq 0;\\
\mu_{1}(u,\chi_{4}, *) & =\textstyle\frac{1}{15}(3 \nu_{3a} + 6 \nu_{5a} +  \nu_{5b} - 4 \nu_{5c} + 6 \nu_{5d} + 3 \nu_{15a} - 2 \nu_{15b} + \alpha_6)\geq 0;\\
\mu_{5}(u,\chi_{4}, *) & =\textstyle\frac{1}{15}(-12 \nu_{3a} - 14 \nu_{5a} - 14 \nu_{5b} - 4 \nu_{5c} - 4 \nu_{5d} - 2 \nu_{15a} - 2 \nu_{15b} + \alpha_7)\geq 0;\\
\mu_{6}(u,\chi_{4}, *) & =\textstyle\frac{1}{15}(-6 \nu_{3a} - 12 \nu_{5a} - 2 \nu_{5b} + 8 \nu_{5c} - 12 \nu_{5d} - 6 \nu_{15a} + 4 \nu_{15b} + \alpha_8)\geq 0;\\
\mu_{0}(u,\chi_{6}, *) & =\textstyle\frac{1}{15}(72 \nu_{3a} - 32 \nu_{5a} - 32 \nu_{5b} + 8 \nu_{5c} + 8 \nu_{5d} - 8 \nu_{15a} - 8 \nu_{15b} + \alpha_9)\geq 0;\\
\mu_{5}(u,\chi_{6}, *) & =\textstyle\frac{1}{15}(-36 \nu_{3a} + 16 \nu_{5a} + 16 \nu_{5b} - 4 \nu_{5c} - 4 \nu_{5d} + 4 \nu_{15a} + 4 \nu_{15b} + \alpha_{10})\geq 0;\\
\mu_{0}(u,\chi_{11}, *) & =\textstyle\frac{1}{15}(-72 \nu_{3a} + 8 \nu_{5a} + 8 \nu_{5b} + 8 \nu_{5c} + 8 \nu_{5d} + 8 \nu_{15a} + 8 \nu_{15b} + \alpha_{11})\geq 0;\\
\mu_{5}(u,\chi_{11}, *) & =\textstyle\frac{1}{15}(36 \nu_{3a} - 4 \nu_{5a} - 4 \nu_{5b} - 4 \nu_{5c} - 4 \nu_{5d} - 4 \nu_{15a} - 4 \nu_{15b} + \alpha_{12})\geq 0;\\
\mu_{1}(u,\chi_{16}, *) & =\textstyle\frac{1}{15}(8 \nu_{3a} -  \nu_{3b} + 9 \nu_{5a} - 11 \nu_{5b} + 4 \nu_{5c} - 6 \nu_{5d} + 3 \nu_{15a} - 2 \nu_{15b} + \alpha_{13})\geq 0;\\
\mu_{6}(u,\chi_{16}, *) & =\textstyle\frac{1}{15}(-16 \nu_{3a} + 2 \nu_{3b} - 18 \nu_{5a} + 22 \nu_{5b} - 8 \nu_{5c} \\
& \hspace{150pt} + 12 \nu_{5d} - 6 \nu_{15a} + 4 \nu_{15b} + \alpha_{14})\geq 0,\\
\end{split}}
\]
where the tuples $(\alpha_1,\dots,\alpha_{14})$ are given in the following table:
$$
\small{
\begin{array}{|c|c|c|c|}\hline
   & \chi(u^5)& \chi(u^3)                  & (\alpha_1,\dots,\alpha_{14})                                   \\ \hline
1  & \chi(3a) & \chi(5a)                   & ( 30, 0,  15, 15, 41, 17, 32, 26, 38, 11, 112, 139, 227, 251 ) \\ \hline  
2  & \chi(3a) & \chi(5b)                   & ( 30, 15, 15, 30, 41, 12, 32, 21, 38, 11, 112, 139, 207, 231 ) \\ \hline  
3  & \chi(3a) & \chi(5c)                   & ( 30, 10, 15, 25, 31, 12, 22, 21, 58, 31, 112, 139, 222, 246 ) \\ \hline  
4  & \chi(3a) & \chi(5d)                   & ( 30, 5,  15, 20, 31, 22, 22, 31, 58, 31, 112, 139, 212, 236 ) \\ \hline  
5  & \chi(3b) & \chi(5a)                   & ( 18, 6,  21, 3,  35, 20, 35, 20, 20, 20, 130, 130, 236, 233 ) \\ \hline  
6  & \chi(3b) & \chi(5b)                   & ( 18, 21, 21, 18, 35, 15, 35, 15, 20, 20, 130, 130, 216, 213 ) \\ \hline  
7  & \chi(3b) & \chi(5c)                   & ( 18, 16, 21, 13, 25, 15, 25, 15, 40, 40, 130, 130, 231, 228 ) \\ \hline  
8  & \chi(3b) & \chi(5d)                   & ( 18, 11, 21, 8,  25, 25, 25, 25, 40, 40, 130, 130, 221, 218 ) \\ \hline  
9  & \chi(3a) & 2\chi(5c)-\chi(5d)         & ( 30, 15, 15, 30, 31, 2,  22, 11, 58, 31, 112, 139, 232, 256 ) \\ \hline  
10 & \chi(3b) & 2\chi(5c)-\chi(5d)         & ( 18, 21, 21, 18, 25, 5,  25, 5,  40, 40, 130, 130, 241, 238 ) \\ \hline  
11 & \chi(3a) & -\chi(5c)+2\chi(5d)        & ( 30, 0,  15, 15, 31, 32, 22, 41, 58, 31, 112, 139, 202, 226 ) \\ \hline  
12 & \chi(3b) & -\chi(5c)+2\chi(5d)        & ( 18, 6,  21, 3,  25, 35, 25, 35, 40, 40, 130, 130, 211, 208 ) \\ \hline  
13 & \chi(3a) & \chi(5a)-\chi(5c)+\chi(5d) & ( 30, -5, 15, 10, 41, 27, 32, 36, 38, 11, 112, 139, 217, 241 ) \\ \hline  
14 & \chi(3b) & \chi(5a)-\chi(5c)+\chi(5d) & ( 18, 1,  21, -2, 35, 30, 35, 30, 20, 20, 130, 130, 226, 223 ) \\ \hline  
15 & \chi(3b) & \chi(5a)+\chi(5b)-\chi(5d) & ( 18, 16, 21, 13, 45, 10, 45, 10, 0,  0,  130, 130, 231, 228 ) \\ \hline  
16 & \chi(3b) & \chi(5a)+\chi(5b)-\chi(5c) & ( 18, 11, 21, 8,  45, 20, 45, 20, 0,  0,  130, 130, 221, 218 ) \\ \hline  
17 & \chi(3a) & \chi(5b)+\chi(5c)-\chi(5d) & ( 30, 20, 15, 35, 41, 2,  32, 11, 38, 11, 112, 139, 217, 241 ) \\ \hline  
18 & \chi(3b) & \chi(5b)+\chi(5c)-\chi(5d) & ( 18, 26, 21, 23, 35, 5,  35, 5,  20, 20, 130, 130, 226, 223 ) \\ \hline  
\end{array}}
$$
In cases 3, 4, 9 and 11 we use the additional inequalities
\[
\begin{split}
\mu_{0}(u,\chi_{2},2) & =\textstyle \frac{1}{15} ( -24 \nu_{3a} + 8 \nu_{5a} + 8 \nu_{5b} - 12 \nu_{5c} - 12 \nu_{5d} - 4 \nu_{15a} - 4 \nu_{15b} - 6)\geq 0;\\
\mu_{0}(u,\chi_{17},7)& =\textstyle \frac{1}{15} ( -64 \nu_{3a} - 16 \nu_{3b} - 8 \nu_{5a} - 8 \nu_{5b} \\
& \hspace{110pt} - 8 \nu_{5c} - 8 \nu_{5d} + 16 \nu_{15a} + 16 \nu_{15b} + 179)\geq 0.\\
\end{split}
\]
In cases 5, 6, 14, 18 we use the additional inequalities
\[
\begin{split}
\mu_{0}(u,\chi_{7},*) & = \textstyle \frac{1}{15} (24 \nu_{3b} + 24 \nu_{5a} + 24 \nu_{5b} - 16 \nu_{5c} - 16 \nu_{5d} + 81) \geq 0;\\
\mu_{0}(u,\chi_{12},7)& = \textstyle \frac{1}{15} (56 \nu_{3a} + 8 \nu_{3b} - 8 \nu_{5a} - 8 \nu_{5b} \\
& \hspace{110pt} - 8 \nu_{5c} - 8 \nu_{5d} + 16 \nu_{15a} + 16 \nu_{15b} + 122) \geq 0.
\end{split}
\]
In cases 7, 8, 10, and 12 we use the additional inequalities
\[
\begin{split}
\mu_{0}(u,\chi_{11},7)&=\textstyle\frac{1}{15}(-56 \nu_{3a} + 16 \nu_{3b} + 8 \nu_{5a} + 8 \nu_{5b} \\
& \hspace{110pt} + 8 \nu_{5c} + 8 \nu_{5d} - 16 \nu_{15a} - 16 \nu_{15b} + 109)\geq 0;\\
\mu_{0}(u,\chi_{12},7)&=\textstyle\frac{1}{15}(56 \nu_{3a} + 8 \nu_{3b} - 8 \nu_{5a} - 8 \nu_{5b} \\
& \hspace{110pt} - 8 \nu_{5c} - 8 \nu_{5d} + 16 \nu_{15a} + 16 \nu_{15b} + 122)\geq 0.\\
\end{split}
\]
In case 13 we use the additional inequalities
\[
\begin{split}
\mu_{3}(u,\chi_{4},*) &=\textstyle\frac{1}{15} (-6 \nu_{3a} - 2 \nu_{5a} - 12 \nu_{5b} - 12 \nu_{5c} + 8 \nu_{5d} + 4 \nu_{15a} - 6 \nu_{15b} + 11) \geq 0; \\
\mu_{6}(u,\chi_{2},2) &=\textstyle\frac{1}{15} (6 \nu_{3a} + 8 \nu_{5a} - 12 \nu_{5b} + 8 \nu_{5c} - 2 \nu_{5d} - 4 \nu_{15a} + 6 \nu_{15b} - 1) \geq 0; \\
\mu_{3}(u,\chi_{7},2) &=\textstyle\frac{1}{15} (16 \nu_{3a} + 4 \nu_{3b} + 12 \nu_{5a} - 28 \nu_{5b} \\
& \hspace{110pt} + 12 \nu_{5c} - 8 \nu_{5d} + 6 \nu_{15a} - 4 \nu_{15b} + 44) \geq 0; \\
\mu_{0}(u,\chi_{17},7)&=\textstyle\frac{1}{15} (-64 \nu_{3a} - 16 \nu_{3b} - 8 \nu_{5a} - 8 \nu_{5b} \\
& \hspace{110pt} - 8 \nu_{5c} - 8 \nu_{5d} + 16 \nu_{15a} + 16 \nu_{15b} + 179) \geq 0.\\
\end{split}
\]
In case 15 we use the additional inequalities
\[
\begin{split}
\mu_{0}(u,\chi_{12}, *) & =\textstyle\frac{1}{15} (128 \nu_{3a} + 8 \nu_{3b} - 40 \nu_{5a} - 40 \nu_{5b} + 8 \nu_{15a} + 8 \nu_{15b} + 122) \geq 0; \\ 
\mu_{3}(u,\chi_{2},2) &=\textstyle\frac{1}{15} (6 \nu_{3a} - 12 \nu_{5a} + 8 \nu_{5b} - 2 \nu_{5c} + 8 \nu_{5d} + 6 \nu_{15a} - 4 \nu_{15b} + 5) \geq 0; \\ 
\mu_{0}(u,\chi_{11},7)&=\textstyle\frac{1}{15} (-56 \nu_{3a} + 16 \nu_{3b} + 8 \nu_{5a} + 8 \nu_{5b} \\
& \hspace{110pt} + 8 \nu_{5c} + 8 \nu_{5d} - 16 \nu_{15a} - 16 \nu_{15b} + 109) \geq 0.\\
\end{split}
\]
In case 16 we use the additional inequalities
\[
\begin{split}
\mu_{0}(u,\chi_{12}, *) & =\textstyle\frac{1}{15}(128 \nu_{3a} + 8 \nu_{3b} - 40 \nu_{5a} - 40 \nu_{5b} + 8 \nu_{15a} + 8 \nu_{15b} + 122)\geq 0;\\
\mu_{0}(u,\chi_{11},7)&=\textstyle\frac{1}{15}(-56 \nu_{3a} + 16 \nu_{3b} + 8 \nu_{5a} + 8 \nu_{5b} \\
& \hspace{110pt} + 8 \nu_{5c} + 8 \nu_{5d} - 16 \nu_{15a} - 16 \nu_{15b} + 109)\geq 0.\\
\end{split}
\]
In case 17 we use the additional inequalities
\[
\begin{split}
\mu_{2}(u,\chi_{2},*) & = \textstyle \frac{1}{15} (5 \nu_{3a} -  \nu_{3b} + 9 \nu_{5a} - 6 \nu_{5b} -  \nu_{5c} + 4 \nu_{5d} - 5) \geq 0; \\ 
\mu_{3}(u,\chi_{2},2) & = \textstyle \frac{1}{15} (6 \nu_{3a} - 12 \nu_{5a} + 8 \nu_{5b} - 2 \nu_{5c} + 8 \nu_{5d} + 6 \nu_{15a} - 4 \nu_{15b} - 1) \geq 0.\\ 
\end{split}
\]
It follows that in all of these cases we have no integral solutions such that all $\mu_{i}(u,\chi_{i},p)$ are non-negative integers.

In case 1 we use the additional inequalities
\[
\begin{split}
\mu_{3}(u,\chi_{4}, *) & =\textstyle\frac{1}{15}(-6 \nu_{3a} - 2 \nu_{5a} - 12 \nu_{5b} - 12 \nu_{5c} + 8 \nu_{5d} + 4 \nu_{15a} - 6 \nu_{15b} + 21) \geq 0; \\ 
\mu_{5}(u,\chi_{12}, *) & =\textstyle\frac{1}{15}(-64 \nu_{3a} - 4 \nu_{3b} + 20 \nu_{5a} + 20 \nu_{5b} - 4 \nu_{15a} - 4 \nu_{15b} + 124) \geq 0; \\ 
\mu_{6}(u,\chi_{2},2)&=\textstyle\frac{1}{15}(6 \nu_{3a} + 8 \nu_{5a} - 12 \nu_{5b} + 8 \nu_{5c} - 2 \nu_{5d} - 4 \nu_{15a} + 6 \nu_{15b} - 6) \geq 0; \\ 
\mu_{3}(u,\chi_{7},2)&=\textstyle\frac{1}{15}(16 \nu_{3a} + 4 \nu_{3b} + 12 \nu_{5a} - 28 \nu_{5b} \\
& \hspace{110pt} + 12 \nu_{5c} - 8 \nu_{5d} + 6 \nu_{15a} - 4 \nu_{15b} + 34) \geq 0; \\ 
\mu_{6}(u,\chi_{7},2)&=\textstyle\frac{1}{15}(16 \nu_{3a} + 4 \nu_{3b} - 28 \nu_{5a} + 12 \nu_{5b} \\
& \hspace{110pt} - 8 \nu_{5c} + 12 \nu_{5d} - 4 \nu_{15a} + 6 \nu_{15b} + 54) \geq 0; \\ 
\mu_{0}(u,\chi_{17},7)&=\textstyle\frac{1}{15}(-64 \nu_{3a} - 16 \nu_{3b} - 8 \nu_{5a} - 8 \nu_{5b} \\
& \hspace{110pt} - 8 \nu_{5c} - 8 \nu_{5d} + 16 \nu_{15a} + 16 \nu_{15b} + 179) \geq 0, \\ 
\end{split}
\]
to obtain only one trivial solution with $\nu_{15b}=1$.

In case 2 we will use the additional inequalities
\[
\begin{split}
\mu_{2}(u,\chi_{2}, *) & =\textstyle\frac{1}{15}(5 \nu_{3a} -  \nu_{3b} + 9 \nu_{5a} - 6 \nu_{5b} -  \nu_{5c} + 4 \nu_{5d} ) \geq 0; \\ 
\mu_{5}(u,\chi_{12}, *) & =\textstyle\frac{1}{15}(-64 \nu_{3a} - 4 \nu_{3b} + 20 \nu_{5a} + 20 \nu_{5b} - 4 \nu_{15a} - 4 \nu_{15b} + 124) \geq 0; \\ 
\mu_{3}(u,\chi_{2},2)&=\textstyle\frac{1}{15}(6 \nu_{3a} - 12 \nu_{5a} + 8 \nu_{5b} - 2 \nu_{5c} + 8 \nu_{5d} + 6 \nu_{15a} - 4 \nu_{15b} - 6) \geq 0; \\ 
\mu_{3}(u,\chi_{7},2)&=\textstyle\frac{1}{15}(16 \nu_{3a} + 4 \nu_{3b} + 12 \nu_{5a} - 28 \nu_{5b} \\
& \hspace{110pt} + 12 \nu_{5c} - 8 \nu_{5d} + 6 \nu_{15a} - 4 \nu_{15b} + 54) \geq 0; \\ 
\mu_{6}(u,\chi_{7},2)&=\textstyle\frac{1}{15}(16 \nu_{3a} + 4 \nu_{3b} - 28 \nu_{5a} + 12 \nu_{5b} \\
& \hspace{110pt} - 8 \nu_{5c} + 12 \nu_{5d} - 4 \nu_{15a} + 6 \nu_{15b} + 34) \geq 0; \\ 
\mu_{0}(u,\chi_{17},7)&=\textstyle\frac{1}{15}(-64 \nu_{3a} - 16 \nu_{3b} - 8 \nu_{5a} - 8 \nu_{5b} \\
& \hspace{110pt} - 8 \nu_{5c} - 8 \nu_{5d} + 16 \nu_{15a} + 16 \nu_{15b} + 179) \geq 0, \\ 
\end{split}
\]
to show that it has the unique trivial solution with $\nu_{15a}=1$.

$\bullet$
It remains to prove part (i) of Theorem \ref{T:2}, showing that there are no
elements of orders 14, 21 and 35 in $V(\mathbb Z G)$.
As in the proof of Theorem \ref{T:1}, below we give the table containing
the data describing the constraints on partial augmentations
$\nu_p$ and $\nu_q$ accordingly to (\ref{E:3})--(\ref{E:5}) for
all of these orders. From this table
part (i) of Theorem \ref{T:2} is derived in the same way as in the
proof of Theorem \ref{T:1}.

\smallskip
\centerline{\small{
\begin{tabular}{|c|c|c|c|c|c|c|c|c|c|c|c|c|c|}
\hline
$|u|$&$p$&$q$&$\xi, \; \tau$&$\xi(C_p)$&$\xi(C_q)$&$l$&$m_1$&$m_p$&$m_q$ \\
\hline
   &   &   &                       &    &   & 0 & 86   & -30   & 0 \\
14 & 2 & 7 & $\xi=(4, 8)_{[*]}$    & -5 & 0 & 2 & 86   & 5     & 0 \\
   &   &   &                       &    &   & 7 & 96   & 30    & 0 \\
\hline
   &   &   & $\xi=(4,7)_{[*]}$     & 3  & 0 & 0 & 90   & 36    & 0 \\
21 & 3 & 7 & $\xi=(4,7)_{[*]}$     & 3  & 0 & 7 & 81   & -18   & 0 \\
   &   &   & $\tau=(18,19)_{[*]}$  & 0  & 2 & 1 & 511  & 0     & 2 \\
\hline
35 & 5 & 7 & $\xi=(2,3)_{[*]}$     & 3  & 0 & 0 & 40   & 72    & 0 \\
   &   &   &                       & 3  & 0 & 7 & 25   & -18   & 0 \\
\hline
\end{tabular}
}}

\section{Proof of Theorem \ref{T:3}}

Let $G$ be the third Janko simple group $J_3$.
It is well known \cite{AFG,GAP} that
$|G|=2^7 \cdot 3^5 \cdot 5 \cdot 17 \cdot 19$ and
$exp(G)=2^3 \cdot 3^2 \cdot 5 \cdot 17 \cdot 19$.

The group $G$ only possesses elements of orders
$2$, $3$, $4$, $5$, $6$, $8$, $9$, $10$, $12$, $15$, $17$ and $19$.
Hence, we shall first investigate normalized units of these orders. By
Proposition \ref{P:2}, the order of each torsion unit divides the
exponent of $G$. So, second we consider normalized units of orders
$18$, $20$, $24$, $30$, $34$, $38$, $45$, $51$, $57$, $85$, $95$ and $323$.
We shall prove that units of all these orders except
$18$, $20$, $24$, $30$ and $45$ do not appear in $V(\mathbb ZG)$. We will
omit cases of units of orders $18$, $20$, $24$, $30$ and $45$ since they are
not products of two disctinct primes, so they do not contribute to Kimmerle's conjecture.

Assume that $u$ is a non-trivial normalized unit and consider each case separately.

\noindent $\bullet$ Let $|u|=2$. Since
there is only one conjugacy class in $G$
consisting of elements or order $2$, this case
immediately follows from Propositions \ref{P:5}
and \ref{P:4}.

\noindent $\bullet$  Let $|u|=3$.  By
(\ref{E:1}) and Proposition \ref{P:4} we have
$\nu_{3a}+\nu_{3b}=1$. By Proposition \ref{P:1}
\[
\begin{split}
\mu_{0}(u,\chi_{2},*) & = \textstyle \frac{1}{3} (-10 \nu_{3a} + 8 \nu_{3b} + 85) \geq 0; \\
\mu_{1}(u,\chi_{2},*) & = \textstyle \frac{1}{3} (5 \nu_{3a} - 4 \nu_{3b} + 85) \geq 0; \\
\mu_{0}(u,\chi_{4},2) & = \textstyle \frac{1}{3} (16 \nu_{3a} - 2 \nu_{3b} + 80) \geq 0, \\
\end{split}
\]
and this system only has the ten solutions listed in part (iii) of Theorem \ref{T:3}.

\noindent $\bullet$  Let $|u|=4$.  By
(\ref{E:1}) and Proposition \ref{P:4} we have
$\nu_{2a}+\nu_{4a}=1$. By Proposition \ref{P:1}
\[
\begin{split}
\mu_{0}(u,\chi_{2},*) & = \textstyle \frac{1}{4} (t_1 + 90) \geq 0; \quad
\mu_{2}(u,\chi_{2},*)  = \textstyle \frac{1}{4} (-t_1 + 90) \geq 0; \\
\mu_{0}(u,\chi_{2},3) & = \textstyle \frac{1}{4} (t_2 + 20) \geq 0; \quad
\mu_{2}(u,\chi_{2},3)  = \textstyle \frac{1}{4} (-t_2 + 20) \geq 0, \\
\end{split}
\]
where $t_1=10 \nu_{2a} + 2 \nu_{4a}$ and $t_2=4 \nu_{2a} - 4 \nu_{4a}$.
Solving this system and applying Proposition \ref{P:6}, only the
three solutions listed in part (iv) of Theorem \ref{T:3} remain.

\noindent $\bullet$  Let $|u|=5$. Then
$\nu_{5a}+\nu_{5b}=1$, and we have the system
\[
\begin{split}
\mu_{1}(u,\chi_{8},2) & = \textstyle \frac{1}{5} (3 \nu_{5a} - 2 \nu_{5b} + 322) \geq 0; \\
\mu_{1}(u,\chi_{2},3) & = \textstyle \frac{1}{5} (-3 \nu_{5a} + 2 \nu_{5b} + 18) \geq 0; \\
\mu_{2}(u,\chi_{2},3) & = \textstyle \frac{1}{5} (2 \nu_{5a} - 3 \nu_{5b} + 18) \geq 0, \\
\end{split}
\]
which only has the eight solutions listed in part (v) of Theorem \ref{T:3}.

\noindent $\bullet$ Let $|u|=8$. Then
$\nu_{2a}+\nu_{4a}+\nu_{8a}=1$ and we need to
consider three cases defined by part (iv) of
Theorem \ref{T:3}. For all of them, put
$t_1 = 20 \nu_{2a} + 4 \nu_{4a} - 4 \nu_{8a}$,
$t_2 = 12 \nu_{2a} + 12 \nu_{4a} - 4 \nu_{8a}$
and $t_3 = 8 \nu_{2a} - 8 \nu_{4a}$.

Case 1. $\chi(u^2) = \chi(4a)$. We have the system of inequalities:
\[
\begin{split}
\mu_{0}(u,\chi_{2},*) & = \textstyle \frac{1}{8} (t_1 + 92) \geq 0; \quad \; 
\mu_{4}(u,\chi_{2},*)   = \textstyle \frac{1}{8} (-t_1 + 92) \geq 0; \\
\mu_{0}(u,\chi_{4},*) & = \textstyle \frac{1}{8} (t_2 + 332) \geq 0; \quad
\mu_{4}(u,\chi_{4},*)   = \textstyle \frac{1}{8} (-t_2 + 332) \geq 0; \\
\mu_{0}(u,\chi_{2},3) & = \textstyle \frac{1}{8} (t_3 + 16) \geq 0; \quad \; 
\mu_{4}(u,\chi_{2},3)   = \textstyle \frac{1}{8} (-t_3 + 16) \geq 0; \\
& \mu_{0}(u,\chi_{4},3)   = \textstyle \frac{1}{8} (16 \nu_{2a} - 8 \nu_{8a} + 88) \geq 0, \\
\end{split}
\]
which only has the following nine solutions satisfying Proposition \ref{P:6} such
that all $\mu_{i}(u,\chi_{j},p)$ are non-negative integers:
\[
\begin{split}
(\nu_{2a},\nu_{4a},\nu_{8a}) \in \{ \;
( 2, 4, -5 ), \; & ( 2, 0, -1 ), \; ( -2, 0, 3 ), \; ( 0, 0, 1 ), \; ( 0, 2, -1 ), \\ &
( 2, 2, -3 ), \; ( -2, -4, 7 ), \; ( -2, -2, 5 ), \; ( 0, -2, 3 ) \; \}.
\end{split}
\]

Case 2. $\chi(u^2) = -2 \chi(2a) + 3 \chi(4a)$. Then we obtain the system
\[
\begin{split}
\mu_{0}(u,\chi_{2},*) & = \textstyle \frac{1}{8} (t_1 + 76) \geq 0; \quad \; 
\mu_{4}(u,\chi_{2},*)  = \textstyle \frac{1}{8} (-t_1 + 76) \geq 0; \\
\mu_{0}(u,\chi_{4},*) & = \textstyle \frac{1}{8} (t_2 + 332) \geq 0; \quad 
\mu_{4}(u,\chi_{4},*)  = \textstyle \frac{1}{8} (-t_2 + 332) \geq 0; \\
\mu_{0}(u,\chi_{2},3) & = \textstyle \frac{1}{8} t_3 \geq 0; \qquad \qquad 
\mu_{4}(u,\chi_{2},3)  = \textstyle \frac{1}{8} (-t_3) \geq 0, \\
\end{split}
\]
which has only three solutions satisfying Proposition \ref{P:6} such
that all $\mu_{i}(u,\chi_{j},p)$ are non-negative integers:
$(\nu_{2a},\nu_{4a},\nu_{8a}) \in \{ \; (0, 0, 1), \; (-2, -2, 5), \; (2, 2, -3) \; \}$.

Case 3. $\chi(u^2) =  2 \chi(2a) -   \chi(4a)$. Put $t_4=28 \nu_{2a} - 4 \nu_{4a} - 4 \nu_{8a}$.
Then the system
\[
\begin{split}
\mu_{0}(u,\chi_{2},*) & = \textstyle \frac{1}{8} (t_1 + 108) \geq 0; \quad 
\mu_{4}(u,\chi_{2},*)   = \textstyle \frac{1}{8} (-t_1 + 108) \geq 0; \\
\mu_{0}(u,\chi_{4},*) & = \textstyle \frac{1}{8} (t_2 + 332) \geq 0; \quad
\mu_{4}(u,\chi_{4},*)   = \textstyle \frac{1}{8} (-t_2 + 332) \geq 0; \\
\mu_{0}(u,\chi_{2},3) & = \textstyle \frac{1}{8} (t_3 + 32) \geq 0; \quad \; 
\mu_{4}(u,\chi_{2},3)   = \textstyle \frac{1}{8} (-t_3 + 32) \geq 0; \\
\mu_{0}(u,\chi_{6},3) & = \textstyle \frac{1}{8} (-t_4 + 116) \geq 0; \quad
\mu_{4}(u,\chi_{6},3)   = \textstyle \frac{1}{8} (t_4 + 116) \geq 0; \\
& \mu_{0}(u,\chi_{4},3)   = \textstyle \frac{1}{8} (16 \nu_{2a} - 8 \nu_{8a} + 104) \geq 0, \\
\end{split}
\]
only has the following nine solutions satisfying Proposition \ref{P:6} such
that all $\mu_{i}(u,\chi_{j},p)$ are non-negative integers:
\[
\begin{split}
( \nu_{2a},\nu_{4a},\nu_{8a}) \in \{ \; &
( 2, 6, -7 ), \; ( 2, 4, -5 ), \; ( 0, 4, -3 ), \; ( 2, 0, -1 ), \; ( -2, 0, 3 ), \; \\ &
( 0, 0, 1 ), \; ( -2, -6, 9 ), \; ( 0, 2, -1 ), \; ( 2, 2, -3 ), \; ( -2, 2, 1 ), \; \\ &
( -2, -4, 7 ), \; ( 0, -4, 5 ), \; ( -2, -2, 5 ), \; ( 0, -2, 3 ), \; ( 2, -2, 1 ) \;
\}.
\end{split}
\]

The union of the solutions of all three cases gives us part (vi) of Theorem \ref{T:3}.

\noindent $\bullet$ Let $|u|=17$.  Then
$\nu_{17a}+\nu_{17b}=1$ and we have the system
\[
\begin{split}
\mu_{1}(u,\chi_{3},19) & = \textstyle \frac{1}{17} (t + 110) \geq 0; \quad 
\mu_{1}(u,\chi_{9},19)   = \textstyle \frac{1}{17} (-t + 706) \geq 0; \\
& \mu_{1}(u,\chi_{2},2) = \textstyle \frac{1}{17} (-10 \nu_{17a} + 7 \nu_{17b} + 78) \geq 0; \\
& \mu_{3}(u,\chi_{2},2) = \textstyle \frac{1}{17} (7 \nu_{17a} - 10 \nu_{17b} + 78) \geq 0, \\
\end{split}
\]
where $t=9 \nu_{17a} - 8 \nu_{17b}$,
with the ten solutions listed in part (vii) of Theorem \ref{T:3}.

\noindent $\bullet$ Let $|u|=19$. Then
$\nu_{19a}+\nu_{19b}=1$, and we have the system
\[
\begin{split}
\mu_{1}(u,\chi_{2},*) & = \textstyle \frac{1}{19} (10 \nu_{19a} - 9 \nu_{19b} + 85) \geq 0; \\
\mu_{2}(u,\chi_{2},*) & = \textstyle \frac{1}{19} (-9 \nu_{19a} + 10 \nu_{19b} + 85) \geq 0; \\
\mu_{1}(u,\chi_{5},2) & = \textstyle \frac{1}{19} (11 \nu_{19a} - 8 \nu_{19b} + 84) \geq 0; \\
\mu_{2}(u,\chi_{5},2) & = \textstyle \frac{1}{19} (-8 \nu_{19a} + 11 \nu_{19b} + 84) \geq 0, \\
\end{split}
\]
with the ten solutions listed in part (viii) of Theorem \ref{T:3}.

$\bullet$ It remains to prove that there are no
elements of orders 34, 38, 51, 57, 85, 95 and 323
in $V(\mathbb Z G)$. We give the table with the
data needed to derive part (i) of Theorem
\ref{T:3} in the same way as in the proofs
of Theorem \ref{T:1} and Theorem \ref{T:2}.

\smallskip
\centerline{\small{
\begin{tabular}{|c|c|c|c|c|c|c|c|c|c|c|c|c|c|}
\hline
$|u|$&$p$&$q$&$\xi, \; \tau$&$\xi(C_p)$&$\xi(C_q)$&$l$&$m_1$&$m_p$&$m_q$ \\
\hline
   &   &   &                       &    &   & 0  & 90   & 80    & 0 \\
34 & 2 & 17&  $\xi=(2)_{[*]}$      &  5 & 0 & 1  & 80   & 5     & 0 \\
   &   &   &                       &    &   & 17 & 80   & -80   & 0 \\
\hline
   &   &   &  $\xi=(7)_{[*]}$      & -10& 0 & 0  & 636  & -180  & 0 \\
38 & 2 & 18&  $\xi=(7)_{[*]}$      & -10& 0 & 19 & 656  & 180   & 0 \\
   &   &   &  $\tau=(13)_{[*]}$    &  15& 0 & 1  & 1600 & 15    & 0 \\
\hline
   &   &   &  $\xi=(9)_{[*]}$      &  6 & 0 & 0  & 828  & 192   & 0 \\
51 & 3 & 17&  $\xi=(9)_{[*]}$      &  6 & 0 & 17 & 810  & -96   & 0 \\
   &   &   &  $\tau=(13)_{[*]}$    & -5 & 0 & 0  & 1605 & -160  & 0 \\
\hline
   &   &   &                       &    &   & 0  & 1605 & -180  & 0 \\
57 & 3 & 19&  $\xi=(13)_{[*]}$     & -5 & 0 & 3  & 1605 &   10  & 0 \\
   &   &   &                       &    &   & 19 & 1620 &   90  & 0 \\
\hline
85 & 5 & 17&  $\xi=(4,8)_{[3]}$    & -2 & 0 & 0  & 400  &  -128 & 0 \\
   &   &   &                       &    &   & 17 & 410  &   32  & 0 \\
\hline
95 & 5 & 19&  $\xi=(4)_{[2]}$      &  0 & 4 & 0  & 152  &    0  & 288 \\
   &   &   &                       &    &   & 19 & 152  &    0  & -72 \\
\hline
323& 17& 19&  $\xi=(6)_{[3]}$      &  0 & 1 & 0  & 171  &    0  & 288 \\
   &   &   &                       &    &   & 19 & 171  &    0  & -18 \\
\hline
\end{tabular}
}}

\section{Some remarks on the group $J_4$}\label{Janko4} 

If $G \cong J_4$, then $|G| = 2^{21} \cdot 3^3
\cdot 5 \cdot 7 \cdot 11^3 \cdot 23 \cdot 29
\cdot 31 \cdot 37 \cdot 43$ and $exp(G) = 2^4
\cdot 3 \cdot 5 \cdot 7 \cdot 11 \cdot 23 \cdot
29 \cdot 31 \cdot 37 \cdot 43$ (see
\cite{AFG,GAP}).
 From Propositions \ref{P:5} and \ref{P:4} it immediately follows that
units of orders $3$, $5$, $23$ and $29$ are
rationally conjugate to a group element.
Nevertheless, further computations are limited by
the fact that $p$-Brauer character tables for
$J_4$ are not known for $p \in
\{2,3,11,23,29,31,43\}$ (see
\verb+http://www.math.rwth-aachen.de/~MOC/work.html+).
For example, for units of order $31$ the best
restriction that can be obtained by applying
Proposition \ref{P:1} to all known ordinary and
$p$-Brauer character tables is given by the
system
\[
\begin{split}
\mu_1(u,\chi_{56},*) & = \textstyle \frac{1}{31}
( 21 \nu_{31a} - 10 \nu_{31b} - 10 \nu_{31c} +
2001151845 ) \geq 0; \\ \mu_3(u,\chi_{56},*) & =
\textstyle \frac{1}{31} ( -10 \nu_{31a} - 10
\nu_{31b} + 21 \nu_{31c} + 2001151845 ) \geq 0;
\\ \mu_5(u,\chi_{56},*) & = \textstyle
\frac{1}{31} ( -10 \nu_{31a} + 21 \nu_{31b} - 10
\nu_{31c} + 2001151845 ) \geq 0, \\ \end{split}
\] where $\nu_{31a}+\nu_{31b}+\nu_{31c}=1$, and
all $\mu_i(u,\chi_{j},*)$ are non-negative
integers for every tuple
$(\nu_{31a},\nu_{31b},\nu_{31c})$ such that
$\nu_{31b}, \nu_{31c} \geq -64553285$ and
$\nu_{31b}+\nu_{31c} \leq 64553286$ (thus, it
will have 18,752,070,203,460,153 solutions which
is too far from rational conjugacy).

Hopefully, further progress can be made
using the Luthar-Passi method if $p$-Brauer
character tables for the missing values of $p$
will become available.

\subsection*{Acknowledgments}
The authors are grateful to Ian Gent, Tom Kelsey and
Andrea Rendl for their advice in using constraint
programming solvers.

\bibliographystyle{amsplain}
\bibliography{BovdiJespersKonovalovJanko}

\end{document}